\documentclass[a4paper,12pt,reqno]{amsart}
\usepackage{amsmath,graphicx,graphics,amssymb}
\newtheorem{con}{Conjecture}
\newtheorem{theo}{Theorem}

\numberwithin{equation}{section}

\newcommand{\M}{\operatorname{M}}

\newcommand{\h}{\operatorname{H}}
\newcommand{\epf}{\hfill{$\square$}\medskip}

\newmuskip\pFqskip
\mathchardef\pFcomma=\mathcode`, 

\begin{document}

\title{Centrally symmetric tilings of fern-cored hexagons}

\author{Mihai Ciucu}
\address{Department of Mathematics, Indiana University, Bloomington, Indiana 47405}

\thanks{Research supported in part by NSF grant DMS-1501052}

\begin{abstract} In this paper we enumerate the centrally symmetric lozenge tilings of a hexagon with a fern removed from its center. The proof is based on a variant of Kuo's graphical condensation method.
An unexpected connection with the total number of tilings is established~---~when suitably normalized, the number of centrally symmetric tilings is equal to the square root of the total number of tilings.
The results we present can be regarded as a new extension of the enumeration of self-complementary plane partitions that fit in a box.
\end{abstract}

\maketitle

\section{Introduction}



This paper extends the results of our recent article \cite{symffd}, in which we presented simple product formulas for the number of centrally symmetric lozenge tilings of a hexagon with a bowtie shape removed from its center.

Our extension concerns the situation when a fern (a chain of aligned triangles of alternating orientation) is removed from the center of the hexagon, to obtain regions called {\it fern-cored hexagons}, which we introduced in \cite{fv}.

An unexpected connection is established with the total number of tilings of such fern-cored hexagons, which we determined in \cite{fv}: It turns out that, when suitably normalized, the number of centrally symmetric tilings of a fern-cored hexagon is equal to the square root of the total number of tilings!

There are three symmetry classes of fern corned hexagons: The base case, when no symmetry is required (which we addressed in \cite{fv}), the vertically symmetric case (worked out by Lai in \cite{Lai}), and the centrally symmetric case. The results of the current paper complete therefore the enumeration of symmetry classes of fern-cored hexagons.

The broader background for our results is the ever widening circle of results in enumerative combinatorics inspired by the work of MacMahon on plane partitions during the late 19th and early 20th centuries, especially his results and conjectures on enumerating plane partitions that fit in a given box, and their symmetry classes (see \cite{MacM}). 

The prototypical result in this circle is MacMahon's original theorem \cite{MacM} on the number of plane partitions that fit in an $x\times y\times z$ box --- well-known (see e.g. \cite{DT}) to be equivalent to enumerating lozenge tilings of a hexagonal region of side-lengths $x$, $y$, $z$, $x$, $y$, $z$ (in cyclic order) on a triangular lattice ---, which states that this number is
\begin{equation}
P(x,y,z)=\prod_{i=1}^x\prod_{j=1}^y\prod_{k=1}^z \frac{i+j+k-1}{i+j+k-2}.
\label{eaa}
\end{equation}
The very remarkable fact that all symmetry classes of lozenge tilings of a hexagon are given by equally beautiful formulas (see \cite{And}\cite{Sta}\cite{Kup}\cite{Ste}\cite{KKZ} and the survey \cite{Bres} for more recent developments) supplies ample motivation for considering the problem of enumerating symmetry classes of generalizations thereof. For a shamrock shape removed from the center of hexagons, this is treated in \cite{ff}\cite{symffa}\cite{symffb}\cite{symffc}\cite{symffd}, and for fern-cored hexagons in \cite{fv}, \cite{Lai} and the present paper.

\begin{figure}[h]
  \centerline{
\hfill
{\includegraphics[width=0.44\textwidth]{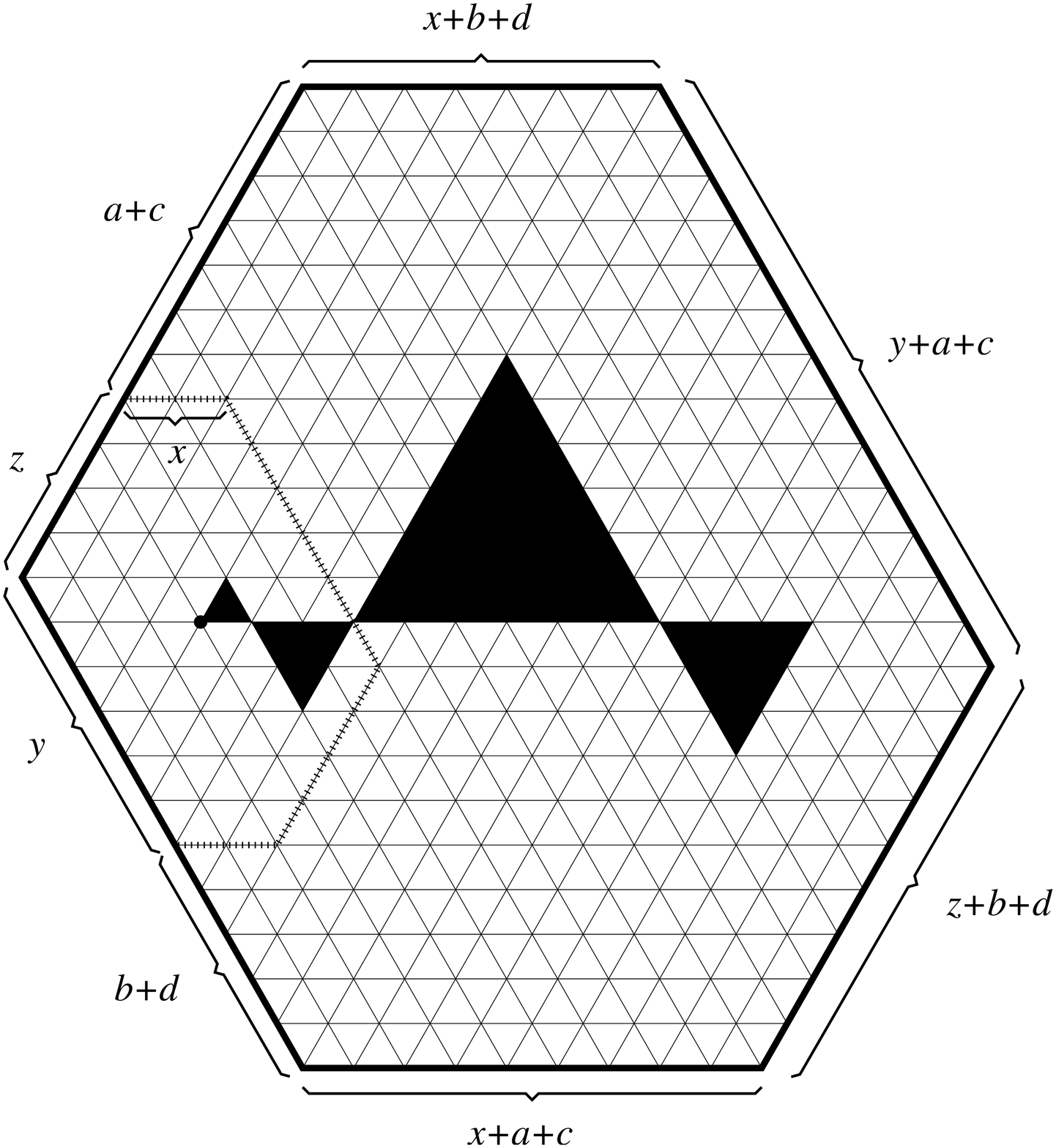}}
\hfill
{\includegraphics[width=0.44\textwidth]{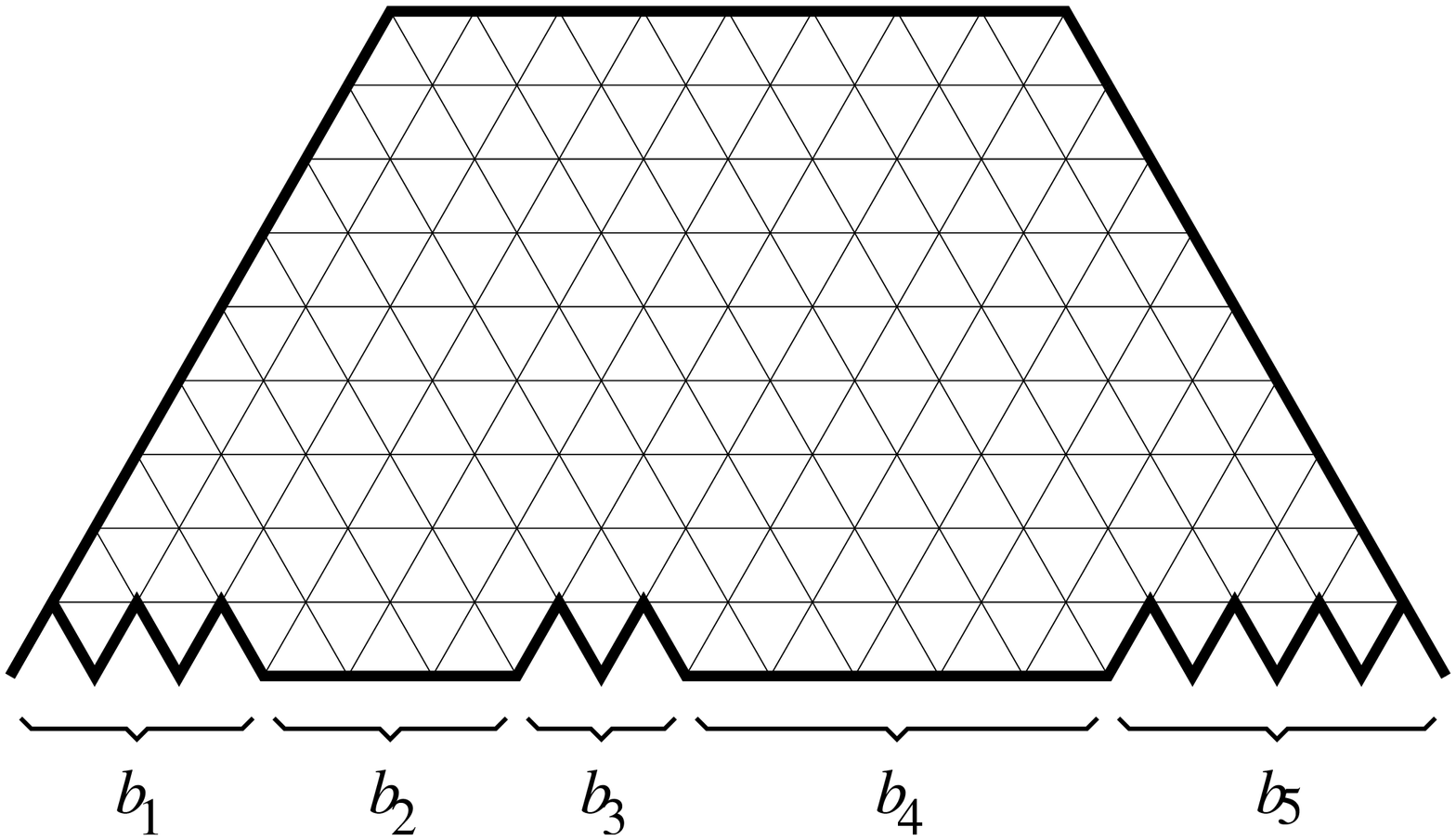}}
}
  \caption{\label{fba} The fern-cored hexagon $FC_{2,6,4}(1,2,6,3)$ (left); the region $S(b_1,b_2,b_3,b_4,b_5)$ for $b_1=3$, $b_2=3$, $b_3=2$, $b_4=5$, $b_5=4$ (right).}
\end{figure}

\section{Statement of main results}

%

Following \cite{fv}, for non-negative integers $a_1,\dotsc,a_k$, we define the {\it fern} $F(a_1,\dotsc,a_k)$ to be a string of $k$ lattice triangles lined up along a horizontal lattice line, touching at their vertices, alternately pointing up and down, having sizes $a_1,\dotsc,a_k$ as encountered from left to right, so that the leftmost triangle is pointing up\footnote{ This does not restrict generality, as the case when the leftmost lobe points down is obtained when $a_1=0$ and $a_2>0$.}; the black structure in Figure \ref{fba} is a 4-lobed fern with lobes of sizes 1, 2, 6 and 2. The {\it base} (or basepoint) of a fern is its leftmost point (the left vertex of the first lobe); for the fern in Figure \ref{fba}, the base is indicated by a black dot.

The fern-cored hexagons are defined as follows. Set
\begin{align}
o&:=a_1+a_3+a_5+\cdots\label{ebaa}\\
e&:=a_2+a_4+a_6+\cdots\label{ebab}
\end{align}
and let $H$ be the hexagon of side-lengths $x+e$, $y+o$, $z+e$, $x+o$, $y+e$, $z+o$ (clockwise from top). We call the hexagon $\bar{H}$ of side-lengths $x$, $y$, $z$, $x$, $y$, $z$ (clockwise from top) that fits in the western corner of $H$ the {\it auxiliary hexagon} (in Figure \ref{fba} the auxiliary hexagon is indicated by a thick dotted line).

Our fern-cored hexagons are obtained by removing a fern from a lattice hexagon such as~$H$. Clearly, such a region can possess centrally symmetric tilings only if the region itself is centrally symmetric. This implies that we must have $o=e$, and also that the center of $H$ is at a lattice point (about which the removed fern is symmetric). Therefore $x$, $y$ and $z$ have the same parity, and the center of $\bar{H}$ is at a lattice point. Then, cf. \cite{fv}, the fern-cored hexagon $FC_{x,y,z}(a_1,\dotsc,a_k)$ is obtained by removing from $H$ the translation of the fern $F(a_1,\dotsc,a_k)$ that places its basepoint at the center of $\bar{H}$ (if the center of $\bar{H}$ is not a lattice point, the fern-cored hexagon is defined slightly differently; the reader can see the details in \cite{fv}).

%
%
%

To state our results, we also need to recall from \cite{fv} that, for non-negative integers $b_1,\dotsc,b_l$, $s(b_1,b_2,\dotsc,b_l)$ denotes the number of lozenge tilings of the semihexagonal region $S(b_1,b_2,\dotsc,b_l)$ with the leftmost $b_1$ up-pointing unit triangles on its base removed, the next segment of length $b_2$ intact, the following $b_3$ removed, and so on (an illustrative example is shown on the right in Figure \ref{fba}; note that the $b_i$'s and the requirement that there are the same number of up-pointing and down-pointing unit triangles in the region --- a necessary condition for the existence of lozenge tilings --- determine the lengths of all four sides of the semihexagon). By the Cohn-Larsen-Propp \cite{CLP} interpretation of the Gelfand-Tsetlin result \cite{GT} we have that
\footnote{ The first equality in \eqref{ebac} holds due to forced lozenges in the tilings of $S(b_1,b_2,\dotsc,b_{2l})$, after whose removal one is left precisely with the region $S(b_1,b_2,\dotsc,b_{2l-1})$.}${}^{,}$\footnote{ We include here the original formula for convenience. Let $T_{m,n}(x_1,\dotsc,x_n)$ be the region obtained from the trapezoid of side lengths $m$, $n$, $m+n$, $n$ (clockwise from top) by removing the up-pointing unit triangles from along its bottom that are in positions $x_1,x_2,\dotsc,x_n$ as counted from left to right. Then
\begin{equation}
\M(T_{m,n}(x_1,\dotsc,x_n))=\prod_{1\leq i<j\leq n}\frac{x_j-x_i}{j-i}.
\label{ebacc}
\end{equation}
}
\begin{align}
&
s(b_1,b_2,\dotsc,b_{2l})=s(b_1,b_2,\dotsc,b_{2l-1})=
\frac{\prod_{\text{\rm $1\leq i\leq j\leq 2l-1$, $j-i+1$ odd}}\h(b_i+b_{i+1}+\cdots+b_j)}{\prod_{\text{\rm $1\leq i\leq j\leq 2l-1$, $j-i+1$ even}}\h(b_i+b_{i+1}+\cdots+b_j)}
\nonumber
\\[5pt]
&\ \ \ \ \ \ \ \ \ \ \ \ \ \ \ \ \ \ \ \ \ \ \ \ \ \ \ \ \ \ \ \ \ \ \ \ \ \ \ \ \ \ \ \ \ \ \ \ 
\times\frac{1}{\h(b_1+b_3+\cdots+b_{2l-1})},
\label{ebac}
\end{align}
where the hyperfactorial $\h(n)$ is defined by
\begin{equation}
\h(n):=0!\,1!\cdots(n-1)!
\label{ebad}
\end{equation}

For a lattice region $R$ on the triangular lattice, we denote by $\M(R)$ the number of lozenge tilings\footnote{ A lozenge is the union of two unit triangles of the lattice that share an edge; a lozenge tiling of a lattice region $R$ is a covering of $R$ by lozenges with no gaps or overlaps.} of $R$. If the region $R$ is centrally symmetric, $\M_\odot(R)$ denotes the number of centrally symmetric lozenge tilings of $R$.

We can now state the main results of our paper.

As pointed out already, a necessary condition for the region $FC_{x,y,z}(a_1,\dotsc,a_k,a_k,\dotsc,a_1)$ to possess any centrally symmetric tilings is that $x$, $y$ and $z$ have the same parity. It is not hard to see that $x$, $y$ and $z$ have to be in fact even in order for this region to have centrally symmetric tilings. The argument in the special case $k=1$ was presented in \cite{symffd}, and it readily extends\footnote{ We note that the region we denoted by $B_{x,y,z,k}$ in \cite{symffd} is the region $FC_{x-k,y-k,z-k}(k,k)$ in the notation of the current paper.} to arbitrary $k$. The picture on the top left in Figure \ref{fca} shows an example of a centrally symmetric fern-cored hexagon.

\begin{theo}
\label{tba}
Let $x$, $y$ and $z$ be even non-negative integers, and let $a_1,\dotsc,a_k$ be non-negative integers.
Then
\begin{align}
\frac
{\M_{\odot}(FC_{x,y,z}(a_1,\dotsc,a_k,a_k,\dotsc,a_1))}
{\M_{\odot}(FC_{x,y,z}(a_1+\cdots+a_k,a_1+\cdots+a_k))}
&
=
\sqrt
{
\frac
{\M(FC_{x,y,z}(a_1,\dotsc,a_k,a_k,\dotsc,a_1))}
{\M(FC_{x,y,z}(a_1+\cdots+a_k,a_1+\cdots+a_k))}
}
\label{ebae}
\\[5pt]
&
=
s(a_1,\dotsc,a_k,a_k,\dotsc,a_1)
\frac
{\h(\frac{x+z}{2}+a_1+\cdots+a_k)}
{\h(\frac{x+y}{2}+a_1+\cdots+a_k)}
\nonumber
\\[10pt]
&\!\!\!\!\!\!\!\!\!\!\!\!\!\!\!\!\!\!\!\!\!\!\!\!\!\!\!\!\!\!\!\!\!\!\!\!\!\!\!\!\!\!\!\!\!\!\!\!\!\!\!\!\!\!\!\!\!\!\!\!\!\!\!\!\!\!\!\!\!
\times
\prod_{1\leq 2i+1\leq 2k} 
\frac
{\h(\frac{x+y}{2}+a_1+\cdots+a_{2i+1})}
{\h(\frac{x+z}{2}+a_1+\cdots+a_{2i+1})}
\prod_{1< 2i< 2k} 
\frac
{\h(\frac{x+z}{2}+a_1+\cdots+a_{2i})}
{\h(\frac{x+y}{2}+a_1+\cdots+a_{2i})},
\label{ebb}
\end{align}
where $a_{k+i}:=a_{k-i+1}$, $i=1,\dotsc,k$.
\end{theo}

Our next result concerns the case when the side-lengths of the outside hexagon do not have the same parity. To capture this case, let\footnote{ We find it convenient to allow $x$ to equal $-1$; since the side-lengths of $FC'_{x,y,z}(a_1,\dotsc,a_k,a_k,\dotsc,a_1)$ involving $x$ contain in fact $x+1$, this guarantees non-negativity of the side-lengths of this region.} $x\geq-1$, $y\geq0$ and $z\geq0$ be integers of the same parity, and consider the hexagon $H$ of side-lengths $x+a+1$, $y+a$, $z+a$, $x+a+1$, $y+a$, $z+a$ (clockwise from top), where $a=a_1+\dotsc+a_k$, and $a_1,\dotsc,a_k$ are non-negative integers. Define $FC'_{x,y,z}(a_1,\dotsc,a_k,a_k,\dotsc,a_1)$ to be the region obtained from $H$ by removing from it the fern $F(a_1,\dotsc,a_n)$ placed so that its base point is at the center\footnote{ For $x=-1$, the auxiliary hexagon $\bar{H}$ gets replaced by the self-crossing contour obtained by starting at the western corner of $H$, going $z$ units in the polar direction $\pi/3$, one unit in the direction $-\pi$, $y$ units in direction $-\pi/3$, $z$ units in the direction $-2\pi/3$, one unit in the direction 1, and then closing up the contour by going $y$ units in direction $2\pi/3$. As $x$, $y$ and $z$ have the same parity, and $y,z\geq0$, it follows that for $x=-1$ we have $x+y,x+z\geq0$. This in turn implies that the center of the described contour is inside $H$, and thus our regions are well defined also for $x=-1$.} of $\bar{H}$, together with the image of this fern across the center of $H$ (these two removed ferns are separated by a unit lattice segment centered at the center of $H$; see the picture on the top left of Figure \ref{fcc} for an example). The following result generalizes\footnote{ The region denoted by $B'_{x,y,z,k}$ in \cite{symffd} is, in the notation of the current paper, the fern-cored hexagon $FC'_{y-k-1,x-k,z-k}(k,k)$; compare the picture on the right in Figure 1 of \cite{symffd} and the region on the top left of Figure \ref{fcc} in the current paper.}  Theorem 5 of \cite{symffd}.

\begin{theo}
\label{tbb}
Let $x\geq-1$, $y\geq0$ and $z\geq0$ be integers of the same the same parity, and let $a_1,\dotsc,a_k$ be non-negative integers.
Then
\begin{align}
&
\frac
{\M_{\odot}(FC'_{x,y,z}(a_1,\dotsc,a_k,a_k,\dotsc,a_1))}
{\M_{\odot}(FC'_{x,y,z}(a_1+\cdots+a_k,a_1+\cdots+a_k))}
%
=
s(a_1,\dotsc,a_{k-1},a_k+1,a_k,\dotsc,a_1)
\nonumber
\\[5pt]
&\ \ 
\times
\frac
{\h(\frac{x+z}{2}+a_1+\cdots+a_k)}
{\h(\frac{x+y}{2}+a_1+\cdots+a_k)}
\prod_{1\leq 2i+1\leq 2k} 
\frac
{\h(\frac{x+y}{2}+a_1+\cdots+a_{2i+1})}
{\h(\frac{x+z}{2}+a_1+\cdots+a_{2i+1})}
\prod_{1< 2i< 2k} 
\frac
{\h(\frac{x+z}{2}+a_1+\cdots+a_{2i})}
{\h(\frac{x+y}{2}+a_1+\cdots+a_{2i})},
\label{ebc}
\end{align}
where $a_{k+1}=a_k+1$, and $a_{k+i}:=a_{k-i+1}$, $i=2,\dotsc,k$.
\end{theo}

\parindent0pt
\textsc{Remark 1.} In \cite{symffd} we provided product formulas for the denominators of the fractions on the left hand side in \eqref{ebae} and \eqref{ebc} (see Theorems 4 and 5 in \cite{symffd}, respectively). Together with the statements of Theorems \ref{tba} and \ref{tbb} above, these provide explicit product formulas for $\M_{\odot}(FC_{x,y,z}(a_1,\dotsc,a_k,a_k,\dotsc,a_1))$ and $\M_{\odot}(FC'_{x,y,z}(a_1,\dotsc,a_k,a_k,\dotsc,a_1))$.

\parindent15pt
\medskip
Unlike Theorem \ref{tba}, Theorem \ref{tbb} does not state a connection between the number of centrally symmetric tilings and the total number of tilings of the region. Given the striking relationship between the two in Theorem \ref{tba} --- the former is just the square root of the latter --- one feels compelled to test whether the same is true in the context of Theorem \ref{tbb}. One difficulty is that in this context, the total number of lozenge tilings is not available from the previous literature.

We work out an explicit conjectural formula in more generality (see Conjecture \ref{tbd}). Our formula implies in particular that the natural counterpart of \eqref{ebae} holds also in the context of Theorem \ref{tbb} (see Remark 2).
We conjecture that this ``square root phenomenon'' holds in fact more generally, for an arbitrary number of removed ferns (see Conjecture \ref{tbd}).



%
%


Our regions depend on non-negative integers $x$, $y$, $z$ (which control the side-lengths of the outer hexagon), $n$ lists of non-negative integers ${\bold a}^{(i)}$, $i=1,\dotsc,n$ (each giving the lobe structure of a removed fern), and non-negative integers $g_1,\dotsc,g_{n-1}$ (which give the separations between the removed ferns).

Define $u$ (resp., $d$) to be the sum of the sizes of all the up-pointing (resp., down-pointing) lobes in the removed ferns. Set
$g=g_1+\dotsc+g_{n-1}$,
and consider the hexagon $H$ of side-lengths $x+d+g$, $y+u$, $z+d$, $x+u+g$, $y+d$, $z+u$ (clockwise from top). If $x$, $y$ and $z$ have the same parity, define
$FC_{x,y,z}^{g_1,\dotsc,g_{n-1}}({\bold a}^{(1)},\dotsc,{\bold a}^{(n)})$
to be the region obtained from $H$ by removing $n$ ferns as follows: The fern $F_1=F({\bold a}^{(1)})$
placed so that its basepoint is
at the center of the auxilliary hexagon $\bar{H}$,
the fern $F({\bold a}^{(2)})$ based $g_1$ units to the right of the rightmost point $F_1$, and so on, ending with removing a fern $F_n=F({\bold a}^{(n)})$ based $g_{n-1}$ units to the right of the rightmost point of $F_{n-1}$.

If $x$, $y$ and $z$ do not have the same parity (and therefore the center of the auxilliary hexagon $\bar{H}$ is the midpoint of a lattice segment), we define our generalized region by modifying the above definition in the same way as when we defined the original fern-cored hexagon $FC_{x,y,z}(a_1,\dotsc,a_k)$ (which is the special case when $n=1$)
in \cite{fv}.
Namely, if $x$ has parity opposite to the parity of $y$ and $z$, translate the system of ferns half a unit to the left from the position described in the previous paragraph. Similarly, if $z$ has parity opposite to the parity of $x$ and $y$, translate the system of ferns by half a unit in the polar direction $-2\pi/3$, and if $y$ has parity opposite to to the parity of $x$ and $z$, translate it it half a unit in the polar direction $2\pi/3$.

We conjecture that the following extension of Theorem 2.1 of \cite{fv} holds.

\begin{con}
\label{tbc}
Let $x$, $y$, $z$ and $g_1,\dotsc,g_{n-1}$ be non-negative integers. Let ${\bold a}^{(i)}=(a_1^{(i)},\dotsc,a_{k_i}^{(i)})$, $i=1,\dotsc,n$ be $n$ lists of non-negative integers. Assume that the lengths $k_1,\dotsc,k_n$ of these lists are even\footnote{ Just like the assumption that the leftmost lobe of a fern points up, this does not restrict generality --- simply set the last lobe size to zero to get ferns with an odd number of lobes.}.

Set $k=k_1+\cdots+k_n$, and let $r_i$ be the distance between the center of the auxilliary hexagon $\bar{H}$ and the rightmost point of the $i$th lobe of the fern system $($in the left-to-right order$)$, in units equal to the side of a unit lattice triangle, for $i=1,\dotsc,k$. Let ${\bold a}_0^{(i)}$ be the list consisting of the single entry $a_1^{(i)}+\dotsc+a_{k_i}^{(i)}$, for $i=1,\dotsc,n$. 
Then
\begin{align}
\nonumber
\\[-7pt]
&
\frac
{\M(FC_{x,y,z}^{g_1,\dotsc,g_{n-1}}({\bold a}^{(1)},\dotsc,{\bold a}^{(n)}))}
{\M(FC_{x,y,z}^{g_1,\dotsc,g_{n-1}}({\bold a_0}^{(1)},\dotsc,{\bold a_0}^{(n)}))}
=
\nonumber
\\[5pt]
&\ \ \ \ \ \ \ \ \ \ \ \ \ \ \ \ \ \ \ \ \, 
s(a_1^{(1)},\dotsc,a_{k_1}^{(1)},g_1+a_1^{(2)},a_2^{(2)},\dotsc,a_{k_2}^{(2)},\dotsc,g_{n-1}+a_1^{(n)},a_2^{(n)},\dotsc,a_{k_n}^{(n)})
\nonumber
\\[5pt]
&\ \ \ \ \ \ \ \ \ \ \ \ \ \ \ \ \
\times
s(a_2^{(1)},\dotsc,a_{k_1}^{(1)},g_1+a_1^{(2)},a_2^{(2)},\dotsc,a_{k_2}^{(2)},\dotsc,g_{n-1}+a_1^{(n)},a_2^{(n)},\dotsc,a_{k_n}^{(n)})
\nonumber
\\[5pt]
&
\times
\prod_{1\leq 2i+1\leq k}
\frac{\h(\lfloor\frac{x+y}{2}\rfloor+r_{2i+1})}{\h(\lfloor\frac{x+z}{2}\rfloor+r_{2i+1})}
\frac{\h(\lceil\frac{x+y}{2}\rceil+r_k-r_{2i+1})}{\h(\lceil\frac{x+z}{2}\rceil+r_k-r_{2i+1})}
\prod_{1<2i<k}
\frac{\h(\lfloor\frac{x+z}{2}\rfloor+r_{2i})}{\h(\lfloor\frac{x+y}{2}\rfloor+r_{2i})}
\frac{\h(\lceil\frac{x+y}{2}\rceil+r_k-r_{2i})}{\h(\lceil\frac{x+z}{2}\rceil+r_k-r_{2i})},
\label{ebg}
\end{align}
%
where $s$ is given by \eqref{ebac}.
\end{con}

\medskip
\parindent0pt
\textsc{Remark 2.} It follows directly from Theorem \ref{tbb} and the special case $n=2$, $g_1=1$ of the above conjectural formula that the couterpart of the ``square root phenomenon'' expressed in Theorem \ref{tba} holds also in the context of Theorem \ref{tbb}, i.e. for two ferns removed symmetrically about the center and separated by a unit lattice segment.

\parindent15pt

\medskip

\parindent0pt
\textsc{Remark 3.} The special case of Conjecture \ref{tbc} when the system of ferns consists of a single fern\footnote{ This special case of Conjecture \ref{tbc} follows by the arguments that proved Theorem 2.1 of \cite{fv}; see the second paragraph of Section 4.} --- the case originally treated in \cite{fv} --- differs from Theorem 2.1 of \cite{fv} only in the choice of the normalizing denominator on the left hand side of the stated equalities: In the latter the normalizing term corresponds to a 2-lobed fern whose up-pointing (resp., down-pointing) lobe has size equal to the sum of all the up-pointing (resp., down-pointing) lobe sizes of the original fern, while in the former the normalizing term corresponds to a single-lobed fern, whose (up-pointing) lobe has size equal to the sum of {\it all lobe sizes} of the original fern.

\parindent15pt  
To see the connection between the two forms, apply the formula from Conjecture \ref{tbc} to the case when the fern system consists of the single fern $F(a_1,\dotsc,a_k)$ (as pointed out in footnote~10, for this special case the statement of Conjecture \ref{tbc} is proved in \cite{fv}). We obtain that
\begin{align}
&
\frac
{\M(FC_{x,y,z}(a_1,\dotsc,a_k))}
{\M(FC_{x,y,z}(a_1+\cdots+a_k))}
=
s(a_1,\dotsc,a_{k-1})s(a_2,\dotsc,a_{k})
\nonumber
\\[5pt]
&\ \ \ \ \ \ \ \ \ \ \ \ \ \ \ \ 
\times
\prod_{1\leq 2i+1\leq k}
\frac{\h(\lfloor\frac{x+y}{2}\rfloor+a_1+\dotsc+a_{2i+1})}{\h(\lfloor\frac{x+z}{2}\rfloor+a_1+\dotsc+a_{2i+1})}
\frac{\h(\lceil\frac{x+y}{2}\rceil+\overline{a_1+\dotsc+a_{2i+1}})}{\h(\lceil\frac{x+z}{2}\rceil+\overline{a_1+\dotsc+a_{2i+1}})}
\nonumber
\\[5pt]
&\ \ \ \ \ \ \ \ \ \ \ \ \ \ \ \ 
\times
\prod_{1<2i<k}
\frac{\h(\lfloor\frac{x+z}{2}\rfloor+a_1+\dotsc+a_{2i+1})}{\h(\lfloor\frac{x+y}{2}\rfloor+a_1+\dotsc+a_{2i+1})}
\frac{\h(\lceil\frac{x+y}{2}\rceil+\overline{a_1+\dotsc+a_{2i}})}{\h(\lceil\frac{x+z}{2}\rceil+\overline{a_1+\dotsc+a_{2i}})},
\label{ebi}
\end{align}
where $\overline{a_1+\dotsc+a_{i}}$ stands for $a_{i+1}+\cdots+a_k$.

The special case of \eqref{ebi} when the fern is the two-lobed fern $F(a_1+a_3+\cdots,a_2+a_4+\cdots)$ becomes
\begin{align}
&
\frac
{\M(FC_{x,y,z}(a_1+a_3+\cdots,a_2+a_4+\cdots)}
{\M(FC_{x,y,z}(a_1+\cdots+a_k))}
=
\nonumber
\\[5pt]
&\ \ \ \ \ \ \ \ \ \ \ \ \ \ \ \ 
\frac{\h(\lfloor\frac{x+y}{2}\rfloor+a_1+a_3+\dotsc)}{\h(\lfloor\frac{x+z}{2}\rfloor+a_1+a_3+\dotsc)}
\frac{\h(\lceil\frac{x+y}{2}\rceil+a_2+a_4+\dotsc)}{\h(\lceil\frac{x+z}{2}\rceil+a_2+a_4+\dotsc)}.
\label{ebj}
\end{align}
Taking now the ratios of the corresponding sides in \eqref{ebi} and \eqref{ebj} yields the formula in Theorem 2.1 of \cite{fv}.

Note that the only difference between \eqref{ebi} and the formula in Theorem 2.1 of \cite{fv} is that the latter has, in addition to the factors on the right hand side of \eqref{ebi}, also the reciprocal of the expression on the right hand side of  \eqref{ebj}. From this point of view, using a single-lobed fern of size equal to the width of the original fern emerges as an even more natural normalization than the two-lobed fern we used in \cite{fv}. 
  

\medskip
\parindent15pt

We conjecture that the square root phenomenon \eqref{ebae} holds in fact more generally, for an arbitrary centrally symmetric system of ferns.

%

\begin{con}
\label{tbd}  
Let $x,y,z$ be non-negative integers of the same parity, and let $g_1,\dotsc,g_n$ be non-negative integers. Let ${\bold a}^{(i)}=(a_1^{(i)},\dotsc,a_{k_i}^{(i)})$, $i=1,\dotsc,n$ be $n$ lists of non-negative integers, and assume that the lengths $k_1,\dotsc,k_n$ of these lists are even\footnote{ See footnote 5.}.

Let ${\bold b}^{(i)}$ be the list of non-negative integers obtained by reversing the order of the elements of the list ${\bold a}^{(i)}$, for $i=1,\dotsc,n$, and consider the region
\begin{equation*}
FC_{x,y,z}^{g_1,\dotsc,g_{n-1},g_n,g_{n-1},\dotsc,g_1}({\bold a}^{(1)},\dotsc,{\bold a}^{(n)},{\bold b}^{(n)},\dotsc{\bold b}^{(1)}).
\end{equation*}
Set $k=2k_1+\cdots+2k_n$, and let $r_i$ be the distance between the leftmost point of the system of the $2n$ removed ferns and the rightmost point of its $i$th lobe $($in the left-to-right order$)$, in units equal to the side of a unit lattice triangle, for $i=1,\dotsc,k$. Let ${\bold a}_0^{(i)}$ be the list consisting of the single entry $a_1^{(i)}+\dotsc+a_{k_i}^{(i)}$, and let ${\bold b}_0^{(i)}=(0,a_1^{(i)}+\dotsc+a_{k_i}^{(i)})$, for $i=1,\dotsc,n$. 
Then
\begin{align}
\nonumber
\\[-7pt]
&
\frac
{\M_\odot(FC_{x,y,z}^{g_1,\dotsc,g_{n-1},g_n,g_{n-1},\dotsc,g_1}({\bold a}^{(1)},\dotsc,{\bold a}^{(n)},{\bold b}^{(n)},\dotsc{\bold b}^{(1)}))}
{\M_\odot(FC_{x,y,z}^{g_1,\dotsc,g_{n-1},g_n,g_{n-1},\dotsc,g_1}({\bold a}_0^{(1)},\dotsc,{\bold a}_0^{(n)},{\bold b}_0^{(n)},\dotsc{\bold b}_0^{(1)}))}
\nonumber
\\[5pt]
&\ \ \ \ \ \ 
=\sqrt{
\frac
{\M(FC_{x,y,z}^{g_1,\dotsc,g_{n-1},g_n,g_{n-1},\dotsc,g_1}({\bold a}^{(1)},\dotsc,{\bold a}^{(n)},{\bold b}^{(n)},\dotsc{\bold b}^{(1)}))}
{\M(FC_{x,y,z}^{g_1,\dotsc,g_{n-1},g_n,g_{n-1},\dotsc,g_1}({\bold a}_0^{(1)},\dotsc,{\bold a}_0^{(n)},{\bold b}_0^{(n)},\dotsc{\bold b}_0^{(1)}))}
}
\\[5pt]
&\ \ \ \ \ \ 
=
s(a_1^{(1)},\dotsc,a_{k_1}^{(1)},g_1+a_1^{(2)},a_2^{(2)},\dotsc,a_{k_2}^{(2)},\dotsc,g_{n-1}+a_1^{(n)},a_2^{(n)},\dotsc,a_{k_n}^{(n)},
\nonumber
\\[5pt]
&\ \ \ \ \ \ \ \ \ \ \ \ \ \ \ \ \ \ \ \ \ 
g_{n}+a_{k_n}^{(n)},a_{k_n-1}^{(n)},\dotsc,a_{1}^{(n)},\dotsc,g_2+a_{k_2}^{(2)},a_{k_2-1}^{(2)},\dotsc,a_{1}^{(2)},g_{1}+a_{k_1}^{(1)},a_{k_1-1}^{(1)},\dotsc,a_{1}^{(1)})
\nonumber
\\[5pt]
&\ \ \ 
\times
\prod_{1\leq 2i+1\leq k}
\frac{\h(\frac{x+y}{2}+r_{2i+1})}{\h(\frac{x+z}{2}+r_{2i+1})}
\prod_{1<2i<k}
\frac{\h(\frac{x+z}{2}+r_{2i})}{\h(\frac{x+y}{2}+r_{2i})}.
\label{ebk}
\end{align}
\end{con}


%





\medskip
Our proof of Theorem \ref{tba} is based on an unusual variant of Kuo's graphical condensation method \cite{KuoOne}\cite{KuoTwo}, which we presented in Theorem 2 of \cite{symffd}. We include it below for convenience. 


\begin{theo} \cite[Theorem 2]{symffd}
\label{kuovar}
Let $G$ be a weighted, centrally symmetric, planar bipartite graph embedded in an annulus. Let $\{a_1,a_2\}$,  $\{b_1,b_2\}$ and $\{c_1,c_2\}$ be pairs of symmetric vertices on the outer face of~$G$. Assume that $a_1$, $b_1$, $c_1$, $a_2$, $b_2$, $c_2$ appear in this cyclic order around the outer face, and that they alternate in color. Let $\{d_1,d_2\}$ be a pair of symmetric vertices of $G$ on the central face. Then we have
\begin{align}
\M_{\odot}(G)\M_{\odot}(G_{abcd})=\M_{\odot}(G_{ab})\M_{\odot}(G_{cd})+
\M_{\odot}(G_{ac})\M_{\odot}(G_{bd})+
\M_{\odot}(G_{ad})\M_{\odot}(G_{bc}),
\label{ebm}
\end{align}
where $G_{abcd}=G\setminus\{a_1,a_2,b_1,b_2,c_1,c_2,d_1,d_2\}$, $G_{ab}=G\setminus\{a_1,a_2,b_1,b_2\}$, etc., and $\M_{\odot}(H)$ denotes the number of centrally symmetric perfect matchings of the graph $H$. 

\end{theo}

The proof of Theorem \ref{tbb} relies on a counterpart of the above result, which allows the vertices $d_1$ and $d_2$ to belong to two different faces, provided they share an edge; it is used in the proof of Theorem 5 in \cite{symffd} (it follows from the special case of Theorem 3 of \cite{symffd} when $F_1$ and $F_2$ share an edge; see also the first paragraph in the proof of Theorem 5 in \cite{symffd}). We state it below for convenient reference.

\begin{theo}
\label{kuovarp}
Let $G$ be a weighted, centrally symmetric, planar bipartite graph embedded in a disk, so that $F_1$ and $F_2$ are two adjacent faces that are each other's image through the central symmetry. Let $\{a_1,a_2\}$,  $\{b_1,b_2\}$ and $\{c_1,c_2\}$ be pairs of symmetric vertices on the outer face of~$G$. Assume that $a_1$, $b_1$, $c_1$, $a_2$, $b_2$, $c_2$ appear in this cyclic order around the outer face, and that they alternate in color. Let $d_1\in F_1$ and $d_2\in F_2$ be images of each other through the central symmetry. Then equality \eqref{ebm} holds.

\end{theo}

\section{Proof of Theorems \ref{tba} and \ref{tbb}}

\medskip
Our proofs of Theorem \ref{tba} (resp., Theorem \ref{tbb}) are by induction on convenient linear quantities in the parameters, and use Theorem \ref{kuovar} (resp., its counterpart mentioned in the previous paragraph) at the induction step.

\medskip

\begin{figure}[h]
  \centerline{
\hfill
{\includegraphics[width=0.36\textwidth]{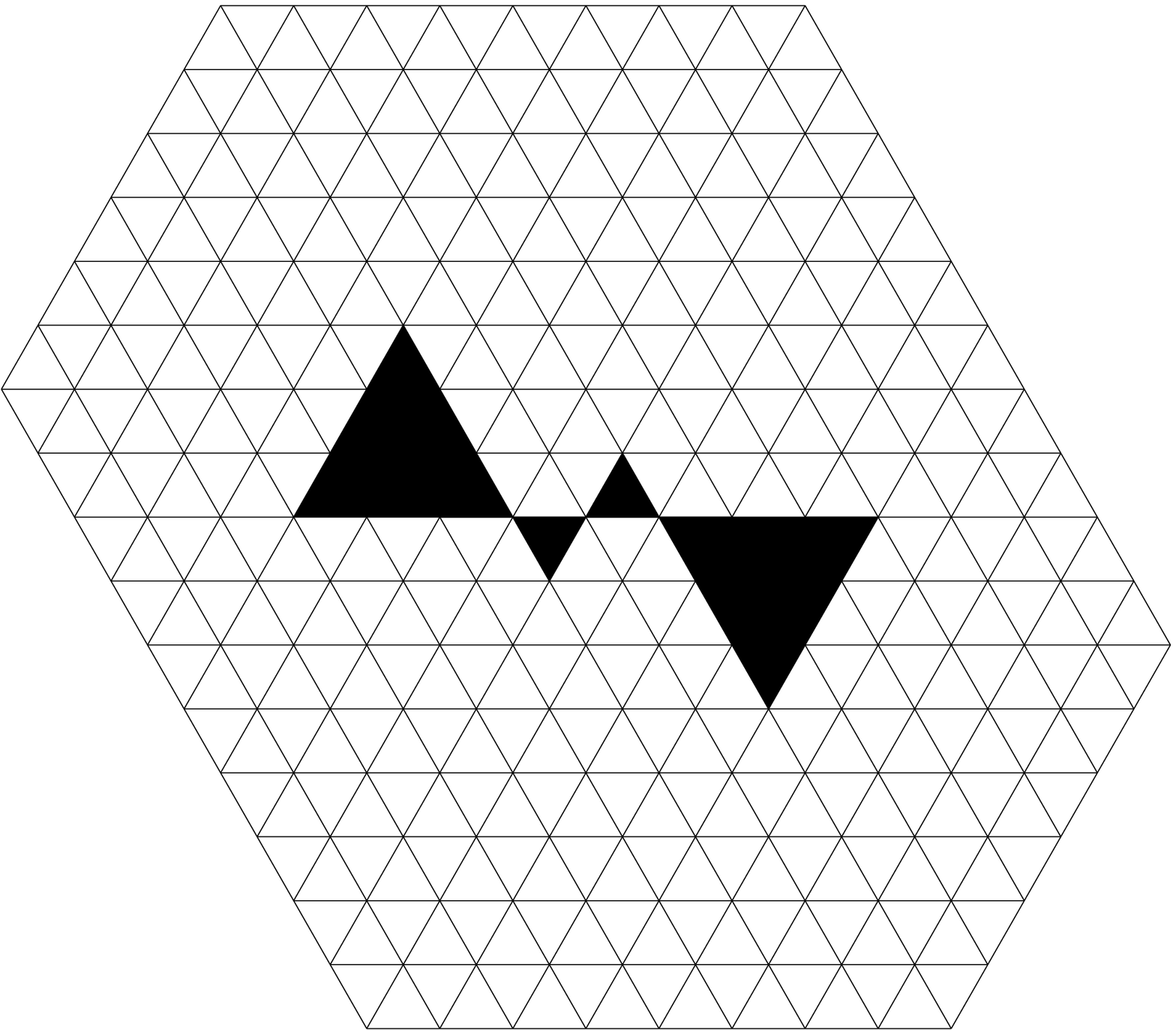}}
\hfill
{\includegraphics[width=0.36\textwidth]{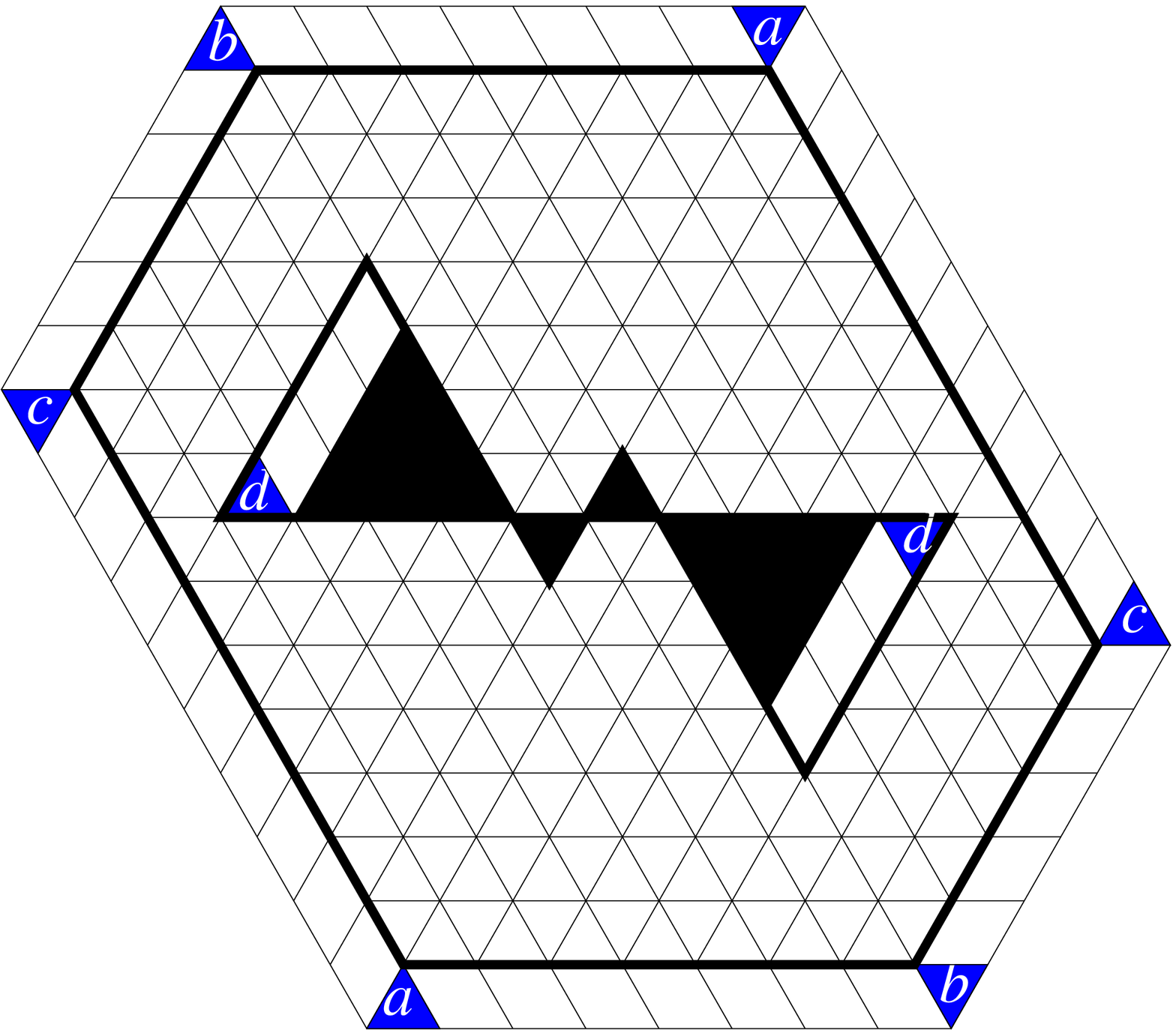}}
\hfill
}
\vskip0.06in
  \centerline{
\hfill
{\includegraphics[width=0.36\textwidth]{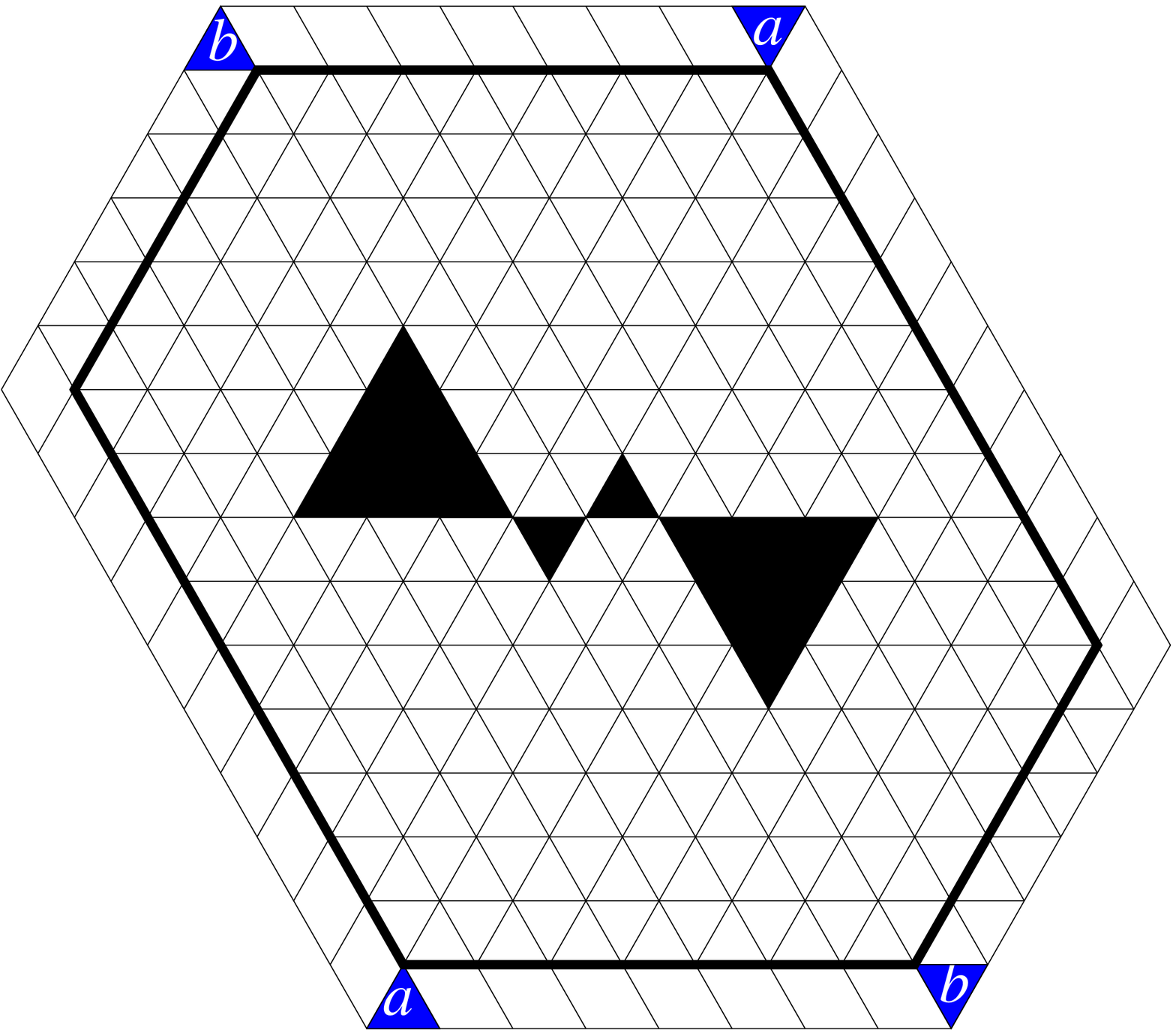}}
\hfill
{\includegraphics[width=0.36\textwidth]{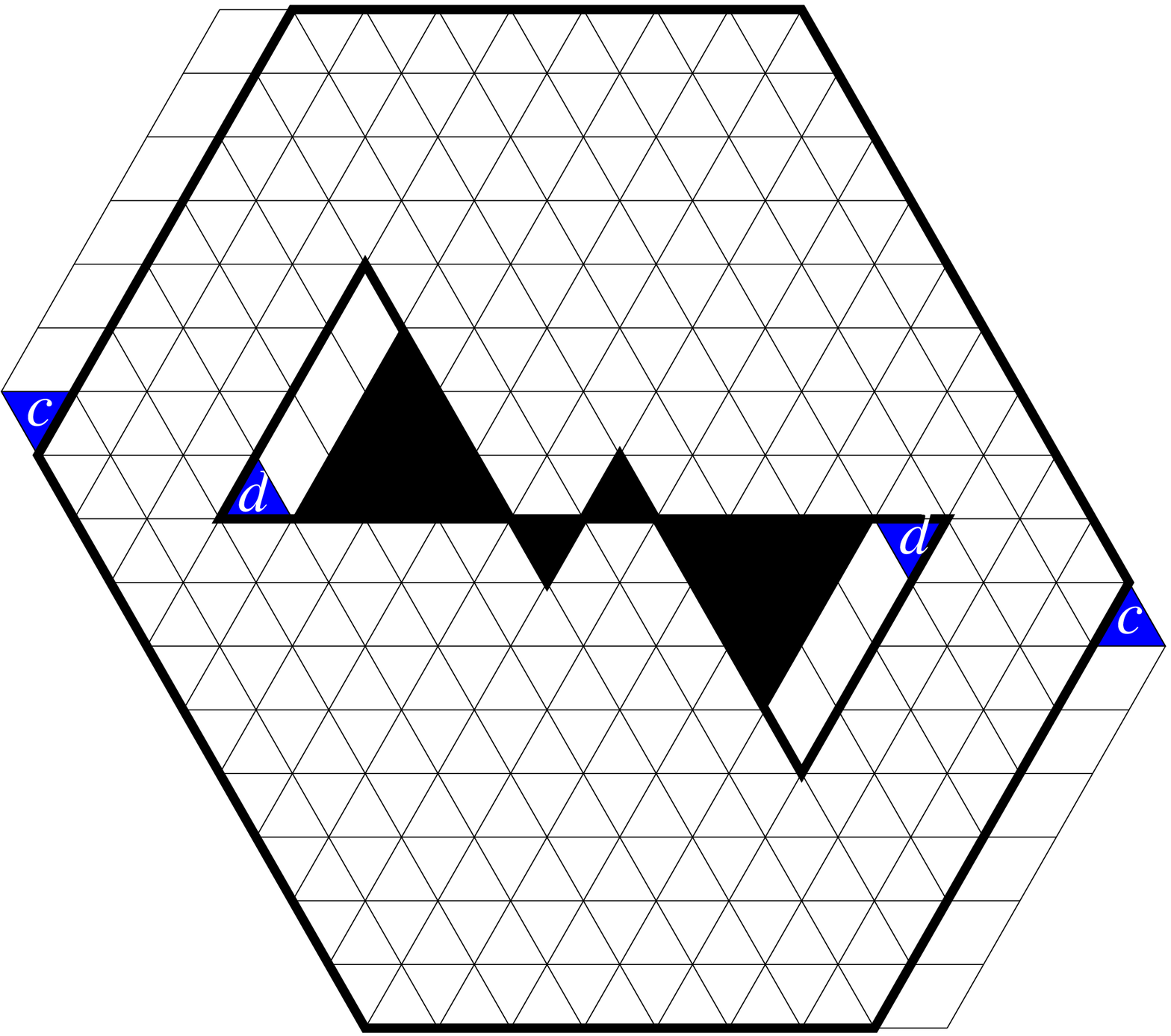}}
\hfill
}
\vskip0.06in
  \centerline{
\hfill
{\includegraphics[width=0.36\textwidth]{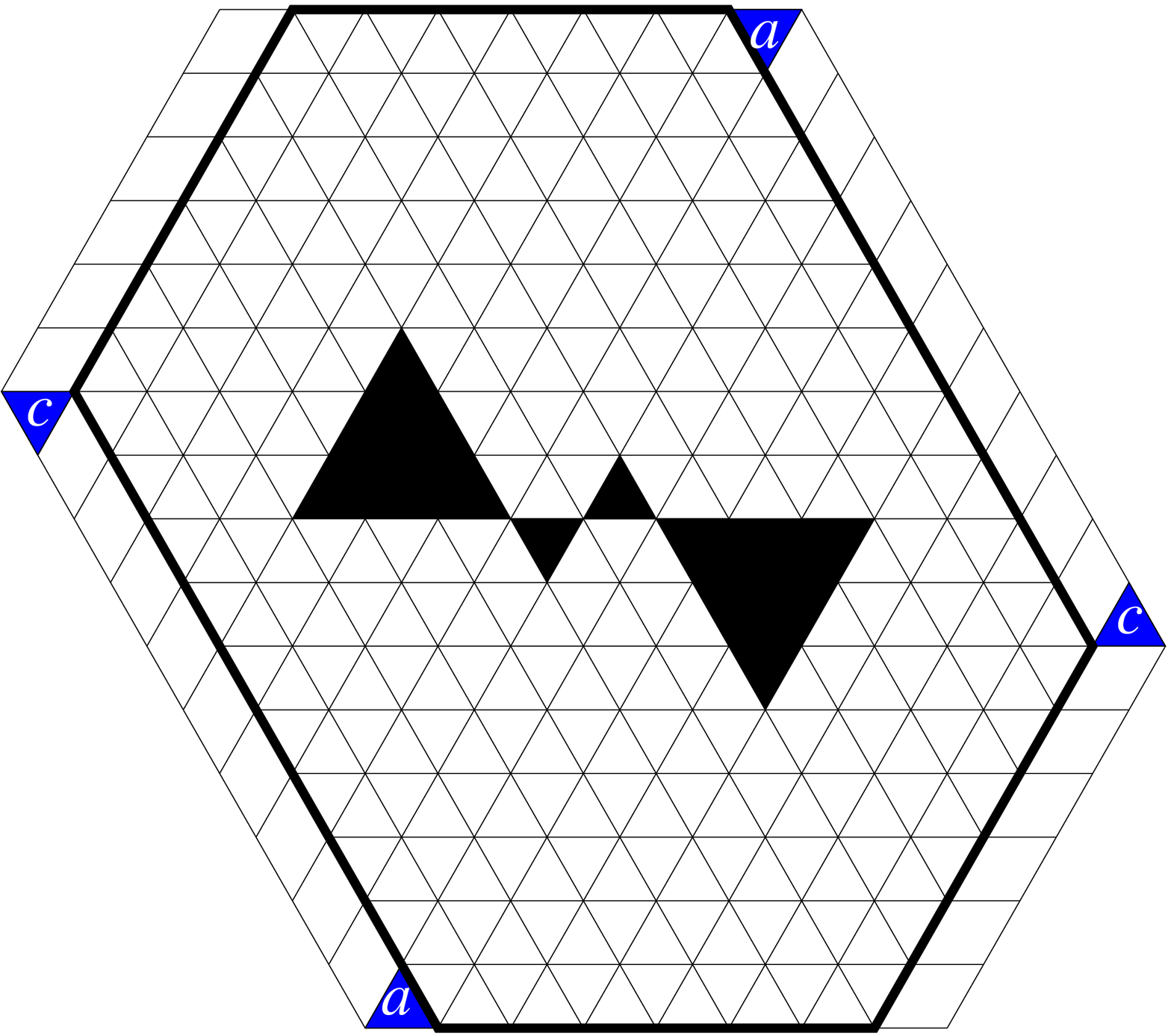}}
\hfill
{\includegraphics[width=0.36\textwidth]{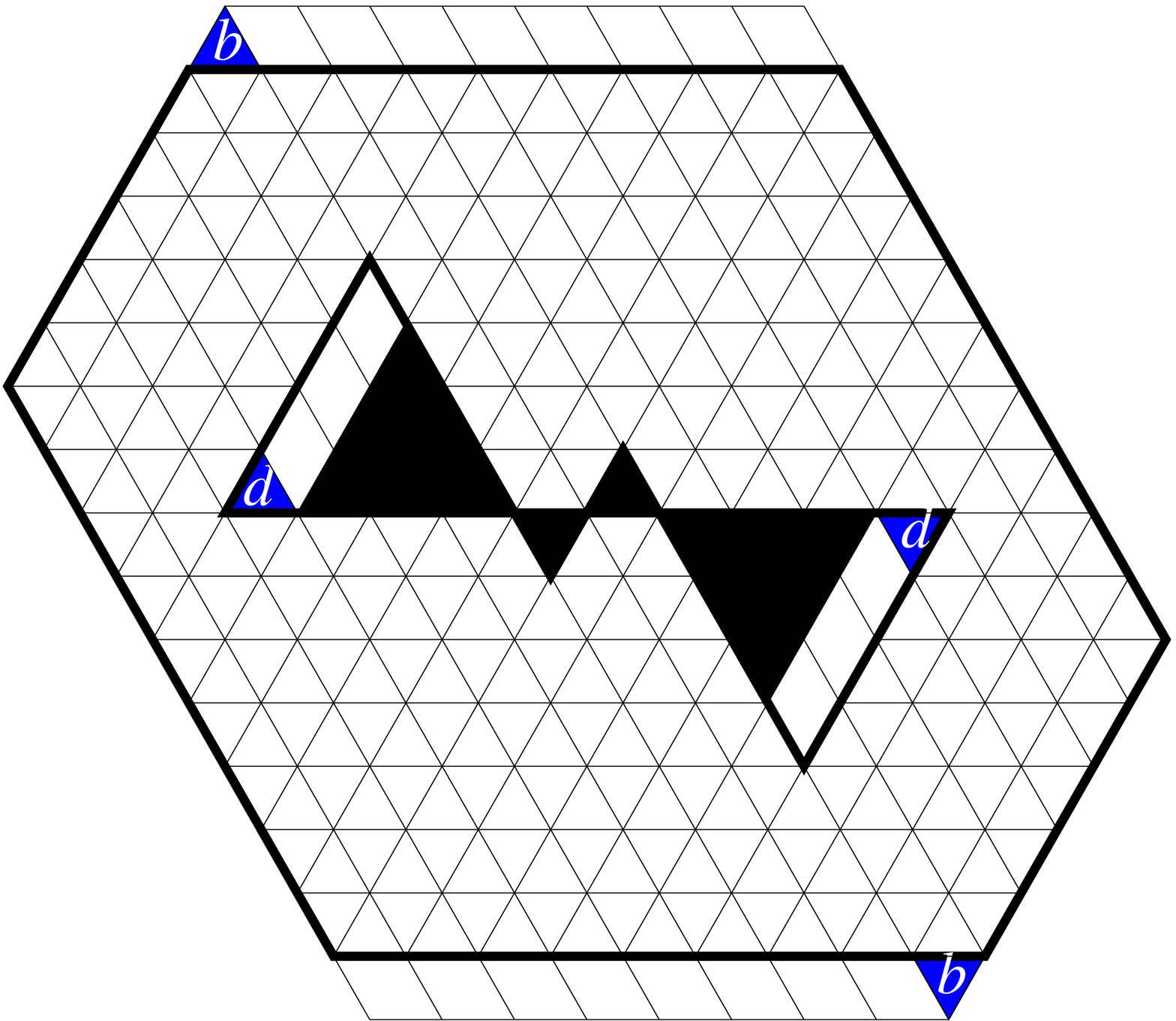}}
\hfill
}
\vskip0.06in
  \centerline{
\hfill
{\includegraphics[width=0.36\textwidth]{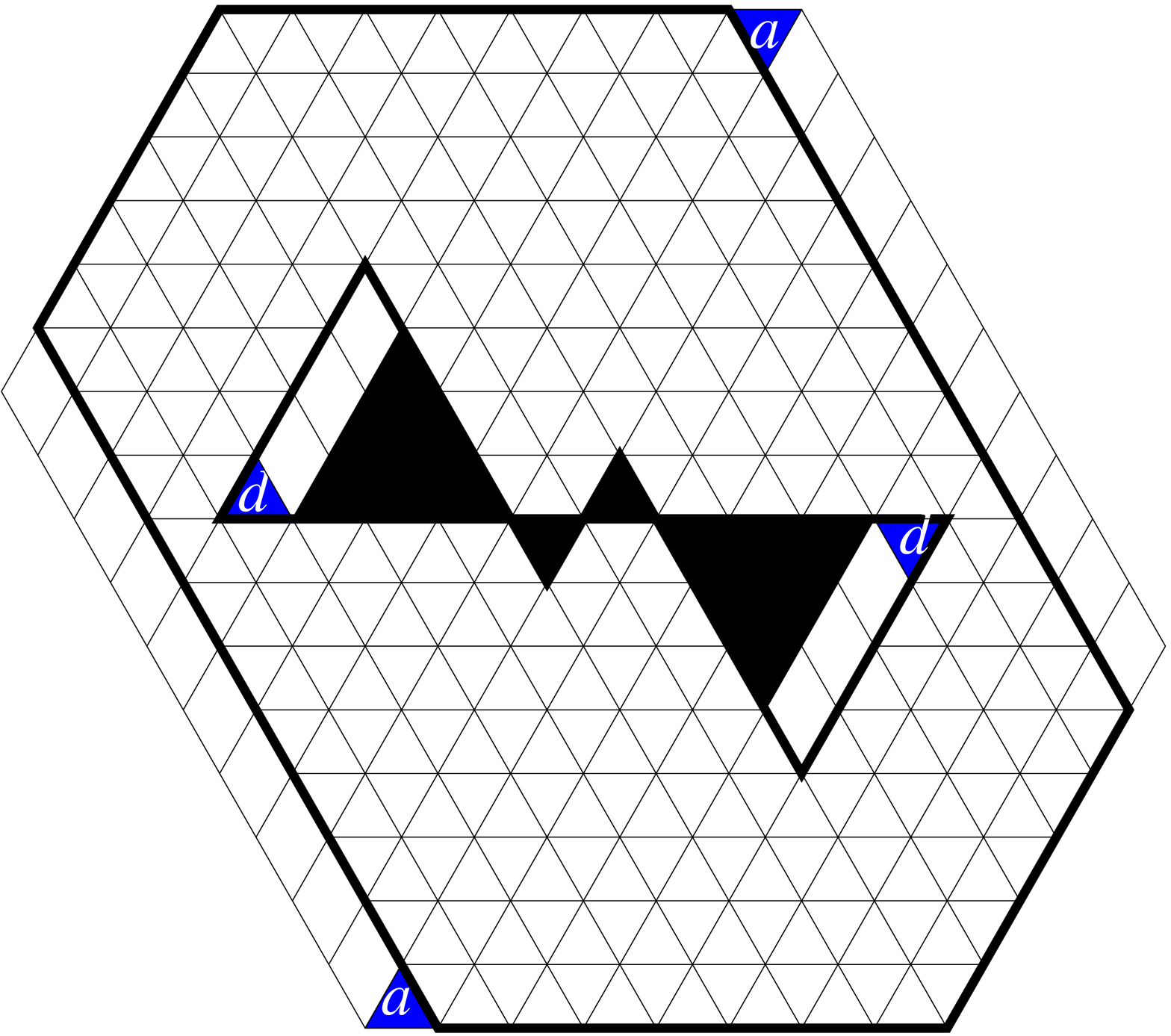}}
\hfill
{\includegraphics[width=0.36\textwidth]{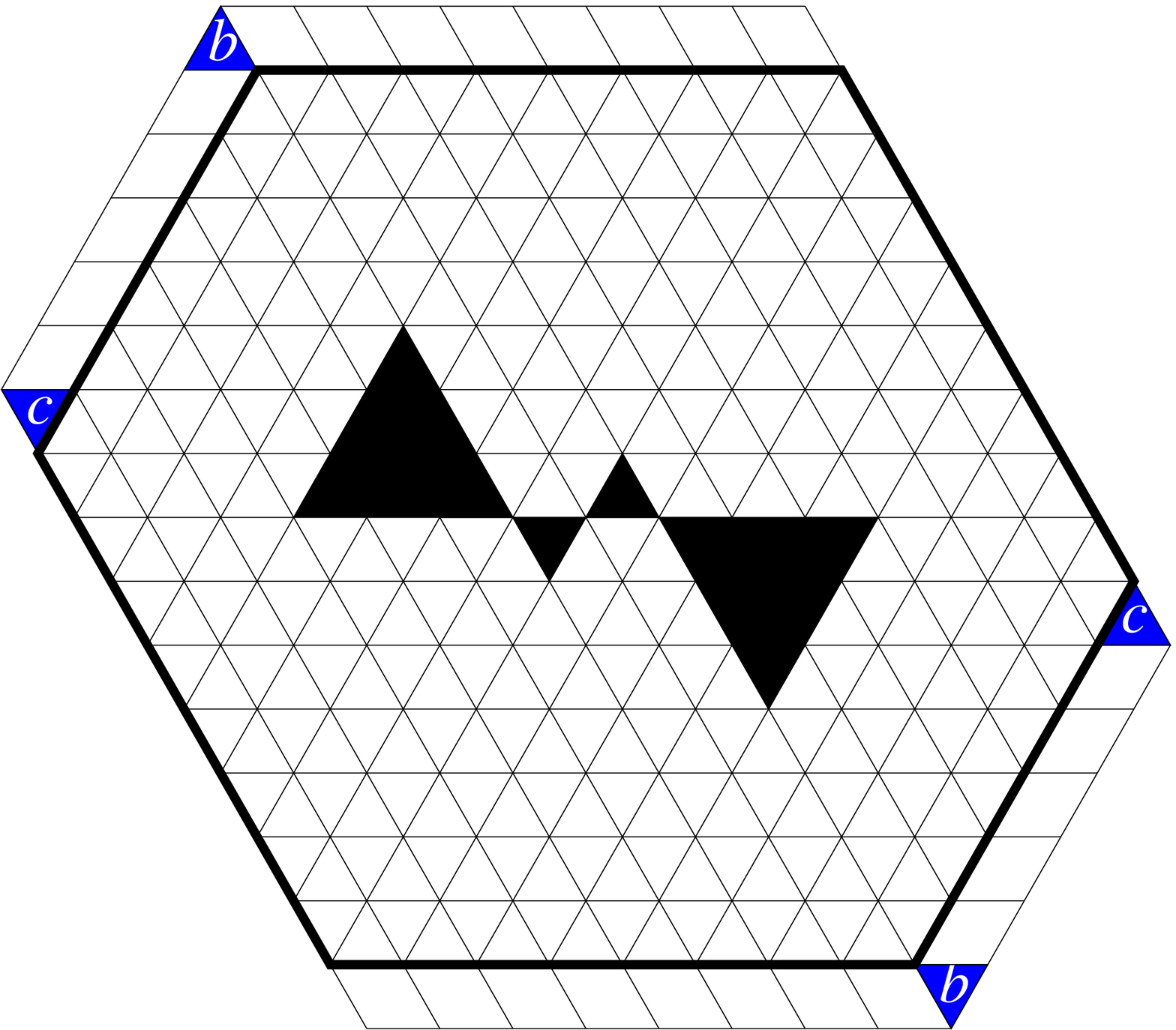}}
\hfill
}
  \caption{\label{fca} Obtaining the recurrence for the number of centrally symmetric tilings of the regions $FC_{x,y,z}(a_1,\dotsc,a_k,a_k,\dotsc,a_1)$.}
\end{figure}

{\it Proof of Theorem \ref{tba}.} Consider the planar dual of the region $FC_{x,y,z}(a_1,\dotsc,a_k,a_k,\dotsc,a_1)$, and apply Theorem \ref{kuovar} to it, with the vertices $a_1,a_2,b_1,b_2,c_1,c_2,d_1,d_2$ of Theorem \ref{kuovar} chosen as indicated on the top right in Figure \ref{fca} ($a_1$ corresponds to the marked unit triangle bottom left, $b_1$, $c_1$, $a_2$, $b_2$, $c_2$ follow in counterclockwise order, and $d_1$ and $d_2$ are next to the fern core). 
After removing the forced lozenges, it follows from Figure \ref{fca} that \eqref{ebm} takes the form
\\[1pt]
\begin{align}
&\!\!
\M_{\odot}(FC_{x,y,z}(a_1,\dotsc,a_k,a_k,\dotsc,a_1))\M_{\odot}(FC_{x-2,y-2,z-2}(a_1+1,a_2\dotsc,a_k,a_k,\dotsc,a_2,a_1+1)
\nonumber
\\[5pt]
&
=\M(FC_{x,y,z-2}(a_1,\dotsc,a_k,a_k,\dotsc,a_1))\M(FC_{x-2,y-2,z}(a_1+1,a_2,\dotsc,a_k,a_k,\dotsc,a_2,a_1+1))
\nonumber
\\[5pt]
&
+\M_{\odot}(FC_{x-2,y,z}(a_1,\dotsc,a_k,a_k,\dotsc,a_1))\M_{\odot}(FC_{x,y-2,z-2}(a_1+1,a_2,\dotsc,a_k,a_k,\dotsc,a_2,a_1+1))
\nonumber
\\[5pt]
&
+
\M_{\odot}(FC_{x-2,y,z-2}(a_1+1,a_2,\dotsc,a_k,a_k,\dotsc,a_2,a_1+1))\M_{\odot}(FC_{x,y-2,z}(a_1,\dotsc,a_k,a_k,\dotsc,a_1)).
\label{eca}
\end{align}
Note that the sum of the $x$-, $y$- and $z$-parameters of the first region in \eqref{eca} is $x+y+z$, and the corresponding sum for the remaining seven regions is stricly less. Equation \eqref{eca} can thus serve to prove the induction step in a proof of \eqref{ebb} by induction on $x+y+z$.

In order for the regions involved in \eqref{eca} to be defined, we need to have $x,y,z\geq2$. The base cases of our induction will thus be the cases when at least one of $x$, $y$ or $z$ is $\leq1$. Since $x$, $y$ and $z$ are even, our base cases are $x=0$, $y=0$ and $z=0$.

\begin{figure}[h]
  \centerline{
\hfill
{\includegraphics[width=0.36\textwidth]{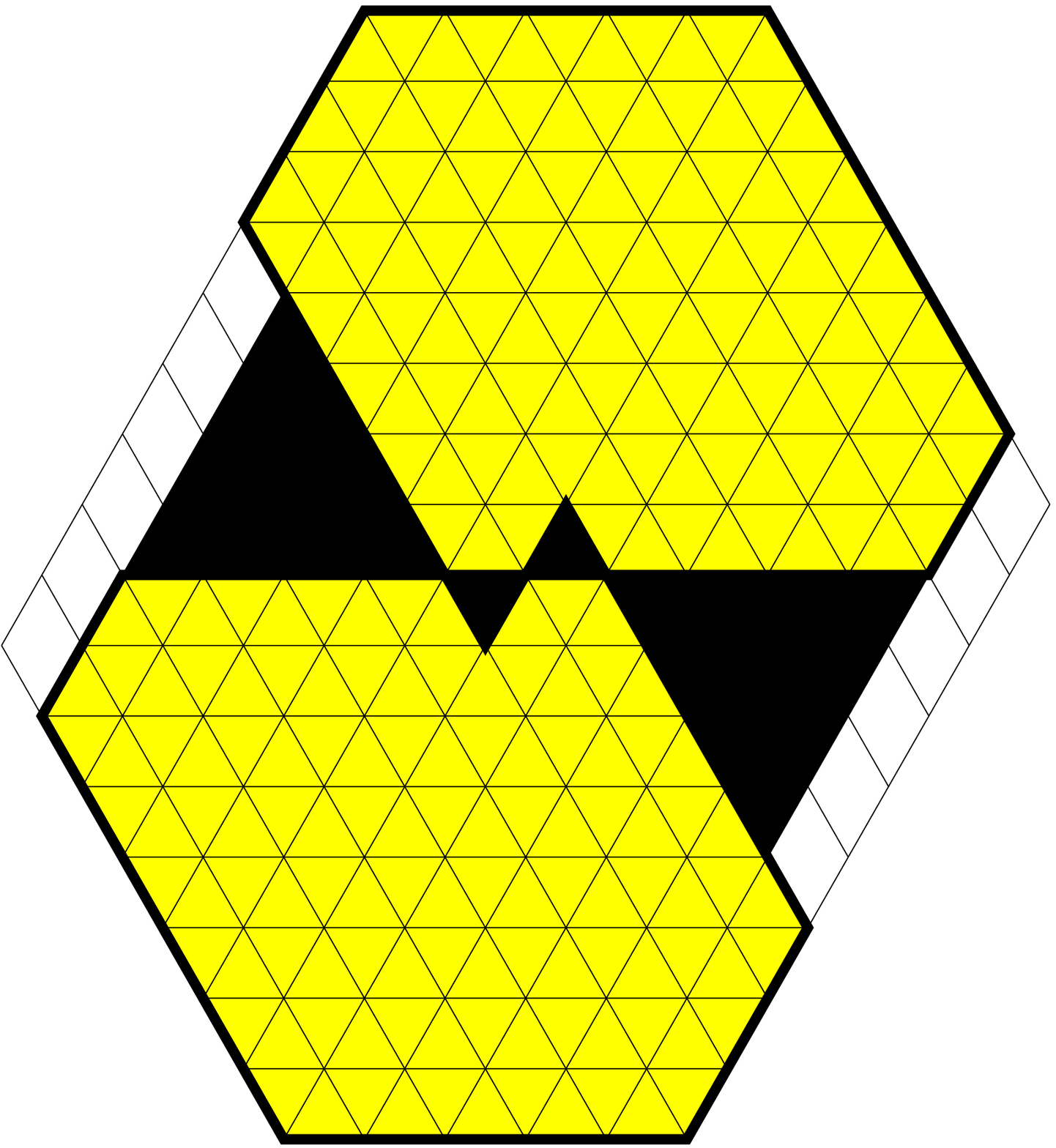}}
\hfill
{\includegraphics[width=0.42\textwidth]{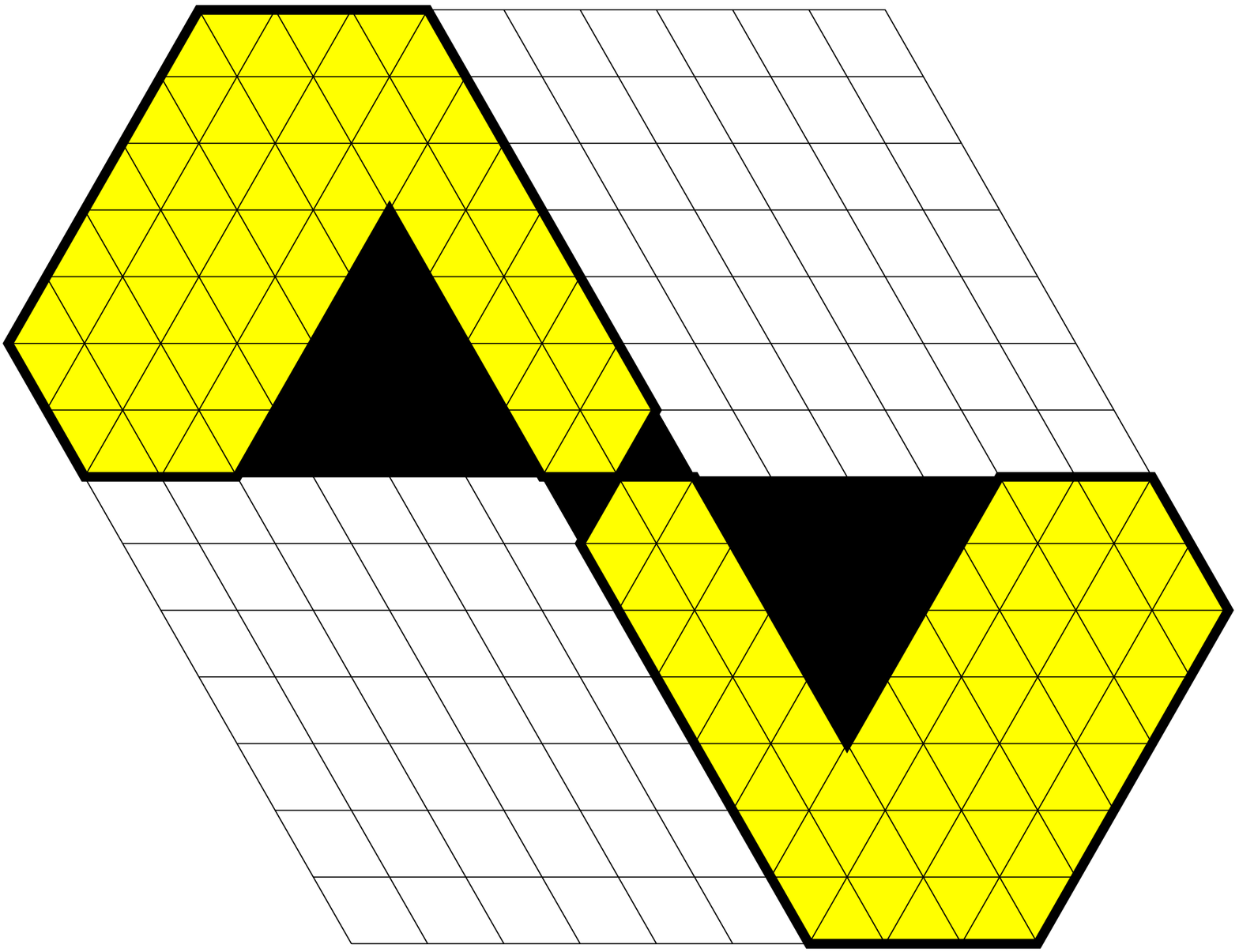}}
\hfill
}
  \caption{\label{fcb} The regions $FC_{0,2,4}(4,1,1,4)$ (left) and $FC_{4,4,0}(4,1,1,4)$ (right).}
\end{figure}

The case $x=0$ is illustrated on the left in Figure \ref{fcb}. In the top shaded region $R$, the length of the top side is equal to the sum of the lengths of the segments that form its base. Encoding lozenge tilings of $FC_{0,y,z}(a_1,\dotsc,a_k,a_k,\dotsc,a_1)$ by families of non-intersecting paths of lozenges that start and end at horizontal lattice segments on its boundary (including the boundary of the removed fern), one readily sees that the previous observation implies that $R$ is internally tiled (i.e., no lozenge crosses its boundary) in each tiling of $FC_{0,y,z}(a_1,\dotsc,a_k,a_k,\dotsc,a_1)$. The same holds for the bottom shaded region, which is the image of $R$ through the symmetry across

In order for the regions involved in \eqref{eca} to be defined, we need to have $x,y,z\geq2$. The base cases of our induction will thus be the cases when at least one of $x$, $y$ or $z$ is $\leq1$. Since $x$, $y$ and $z$ are even, our base cases are $x=0$, $y=0$ and $z=0$.

The case $x=0$ is illustrated on the left in Figure \ref{fcb}. In the top shaded region $R$, the length of the top side is equal to the sum of the lengths of the segments that form its base. Encoding lozenge tilings of $FC_{0,y,z}(a_1,\dotsc,a_k,a_k,\dotsc,a_1)$ by families of non-intersecting paths of lozenges that start and end at horizontal lattice segments on its boundary (including the boundary of the removed fern), one readily sees that the previous observation implies that $R$ is internally tiled (i.e., no lozenge crosses its boundary) in each tiling of $FC_{0,y,z}(a_1,\dotsc,a_k,a_k,\dotsc,a_1)$. The same holds for the bottom shaded region, which is the image of $R$ through the symmetry across the center. Furthermore, tilings of the portion of $FC_{0,y,z}(a_1,\dotsc,a_k,a_k,\dotsc,a_1)$ not covered by these two shaded regions are forced to be as indicated in Figure \ref{fcb}. This implies that
\parindent15pt
\begin{equation*}
\M_\odot(FC_{0,y,z}(a_1,\dotsc,a_k,a_k,\dotsc,a_1))=\M(R).
\end{equation*}  

An explicit product formula for the number of tilings $\M(R)$ of the region $R$ is provided by equation \eqref{ebacc}. Taking the ratio between this product formula and its specialization to the case when the removed fern consists of just two lobes, both of size $a_1+\cdots+a_k$, leads to an explicit product expression for the $x=0$ specialization of the ratio on the left hand side of \eqref{ebb}. One readily sees that this agrees with the $x=0$ specialization of the expression on the right hand side of \eqref{ebb}, which verifies this base case.

The case $y=0$ follows from the above case by symmetry. The case $z=0$ is different, but follows by the same arguments. Indeed, for $z=0$, our region looks as pictured on the right in Figure \ref{fcb}. In the shaded region $R'$ on top, the length of the northwestern side is equal to the sum of the lengths of the portions of its boundary that face southeast. By the argument we used above involving non-intersecting paths of lozenges, it follows that the region $R'$ is always internally tiled. The same is true of the shaded region on the bottom, and as in the case $x=0$ we obtain that
\begin{equation*}
\M_\odot(FC_{x,y,0}(a_1,\dotsc,a_k,a_k,\dotsc,a_1))=\M(R').
\end{equation*}  
Again, we obtain an explicit product formula for $\M(R')$ using equation \eqref{ebacc}, and this leads to a formula for the $z=0$ specialization of the expression on the left hand side of \eqref{ebb}, which is readily verified to agree with the $z=0$ specialization of the expression on the right hand side of \eqref{ebb}. This finishes the verification of the base cases of our induction.

For the induction step, assume that formula \eqref{ebb} holds for all the regions in which the sum of the $x$-, $y$- and $z$-parameters is strictly less than $x+y+z$. We need to deduce that \eqref{ebb} holds in the form it is written.

To see this, use \eqref{eca} to express $\M(FC_{x,y,z}(a_1,\dotsc,a_k,a_k,\dotsc,a_1))$ in terms of the number of lozenge tilings of the other seven regions involved. Since in each of the latter the sum of the $x$-, $y$- and $z$-parameters is strictly less than $x+y+z$, by the induction hypothesis, the number of its tilings is given by the explicit product provided by formula \eqref{ebb}. It is a straightforward calculation to verify that the resulting expression for $\M(FC_{x,y,z}(a_1,\dotsc,a_k,a_k,\dotsc,a_1))$ agrees with the one on the right hand side of \eqref{eca}. Indeed, this is equivalent to the statement that the expression on the right hand side of \eqref{ebb} also satisfies recurrence \eqref{eca}. The details of the verification of an analogous statement are presented in Section 4 of \cite{fv}. As it is shown there, the equality to be verified reduced (due to the fact that certain multipliers in the equation were equal) to the statement of Kuo's graphical condensation theorem (Theorem 2.1 of \cite{KuoOne}) for a special fern-cored hexagon. As it happens (and as it is clear from a comparison of the expression on the right hand sides of \eqref{ebb} and formula (2.4) of \cite{fv}), the multipliers we get in the current set-up are precisely the square roots of the multipliers in Section 4 of \cite{fv}. Therefore, the former are also equal. Using this, one sees that verification of the fact that the expression on right hand side of \eqref{ebb} satisfies recurrence~\eqref{eca} amounts to the special case of\eqref{eca} when the fern-cored hexagon has two lobes of size $a_1+\cdots+a_k$. This completes the proof of formula \eqref{ebb}.

Formula \eqref{ebae} follows now directly from \eqref{ebb} and Theorem 2.1 of \cite{fv}.
\epf

\begin{figure}[h]
  \centerline{
\hfill
{\includegraphics[width=0.36\textwidth]{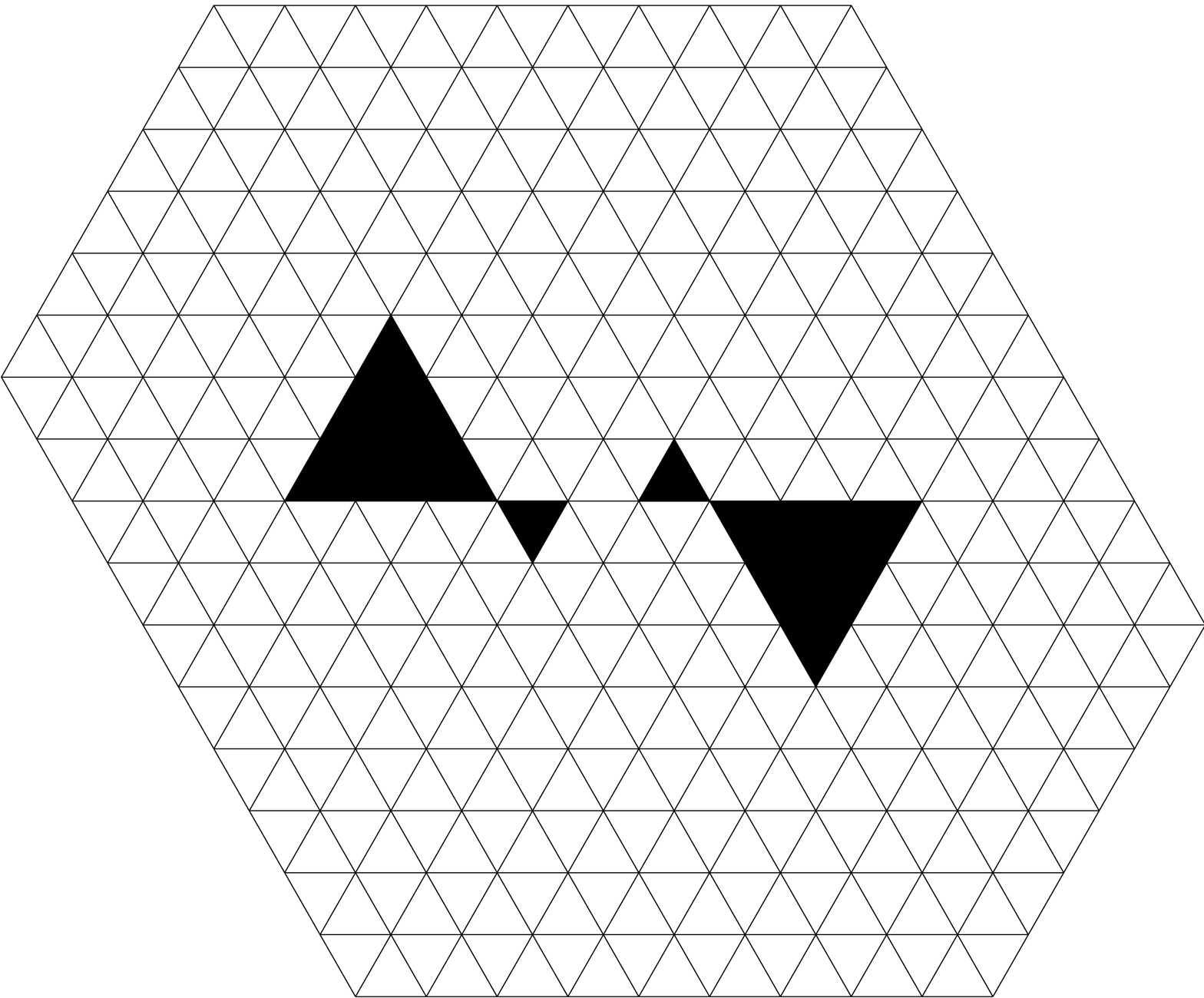}}
\hfill
{\includegraphics[width=0.36\textwidth]{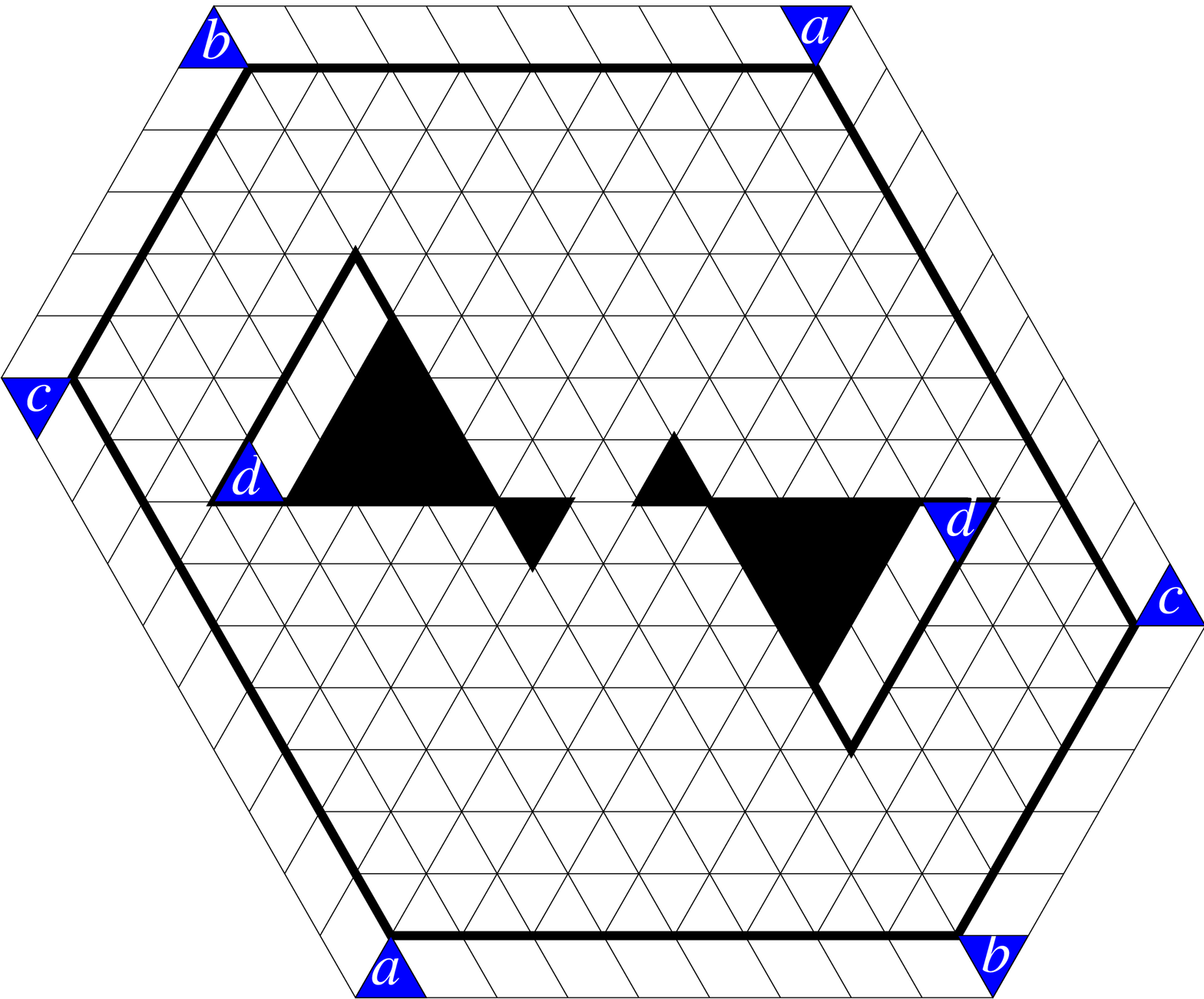}}
\hfill
}
\vskip0.06in
  \centerline{
\hfill
{\includegraphics[width=0.36\textwidth]{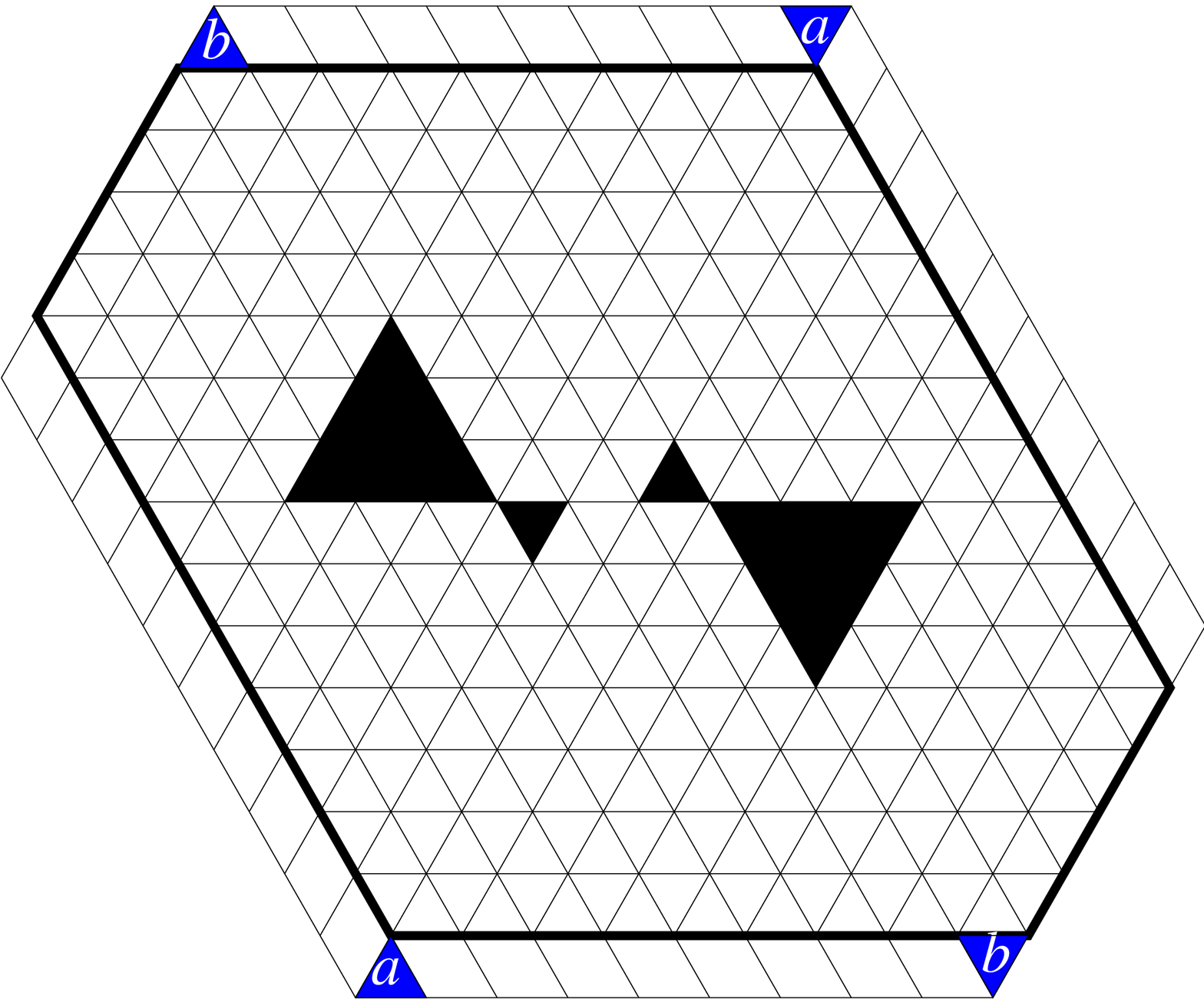}}
\hfill
{\includegraphics[width=0.36\textwidth]{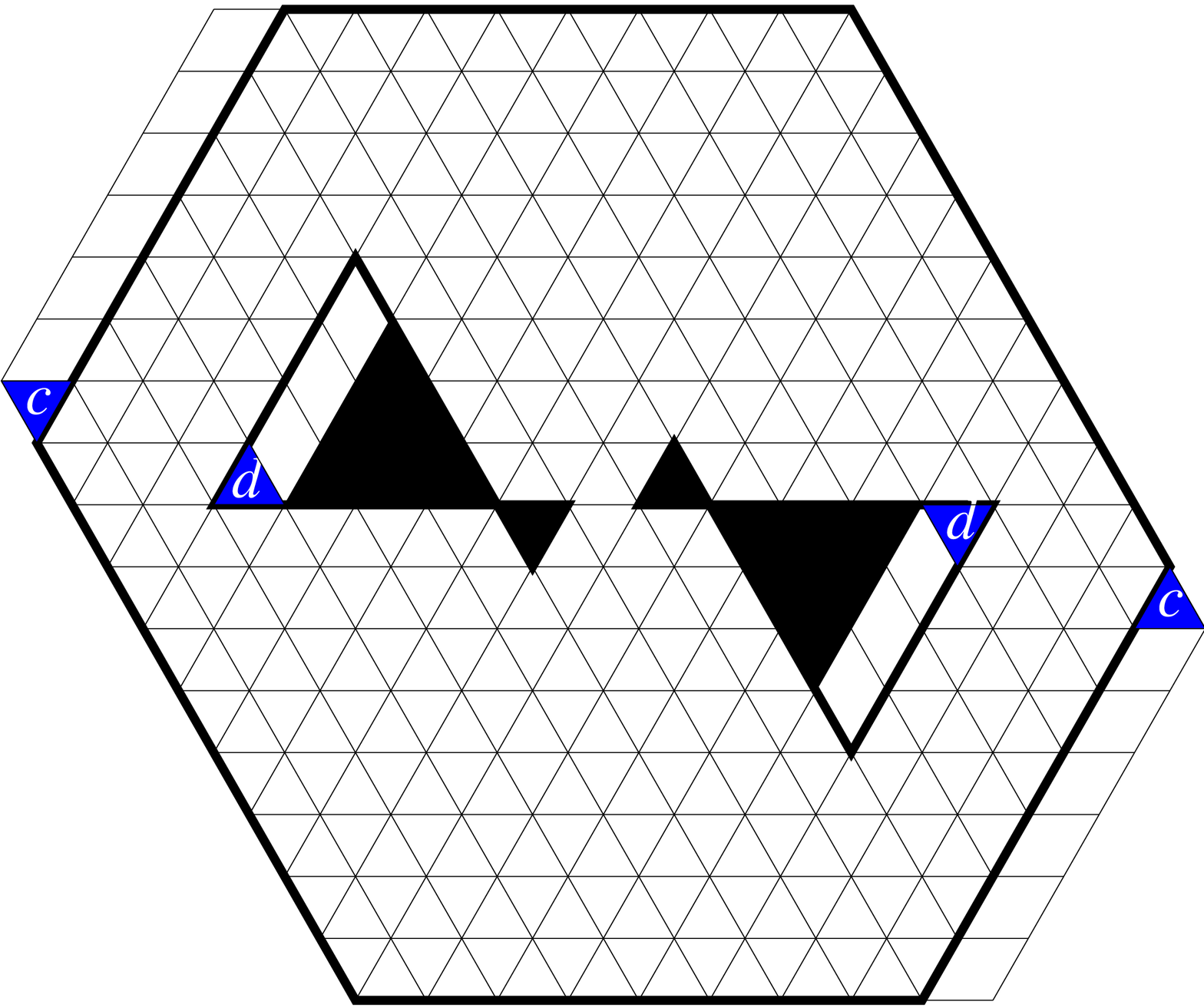}}
\hfill
}
\vskip0.06in
  \centerline{
\hfill
{\includegraphics[width=0.36\textwidth]{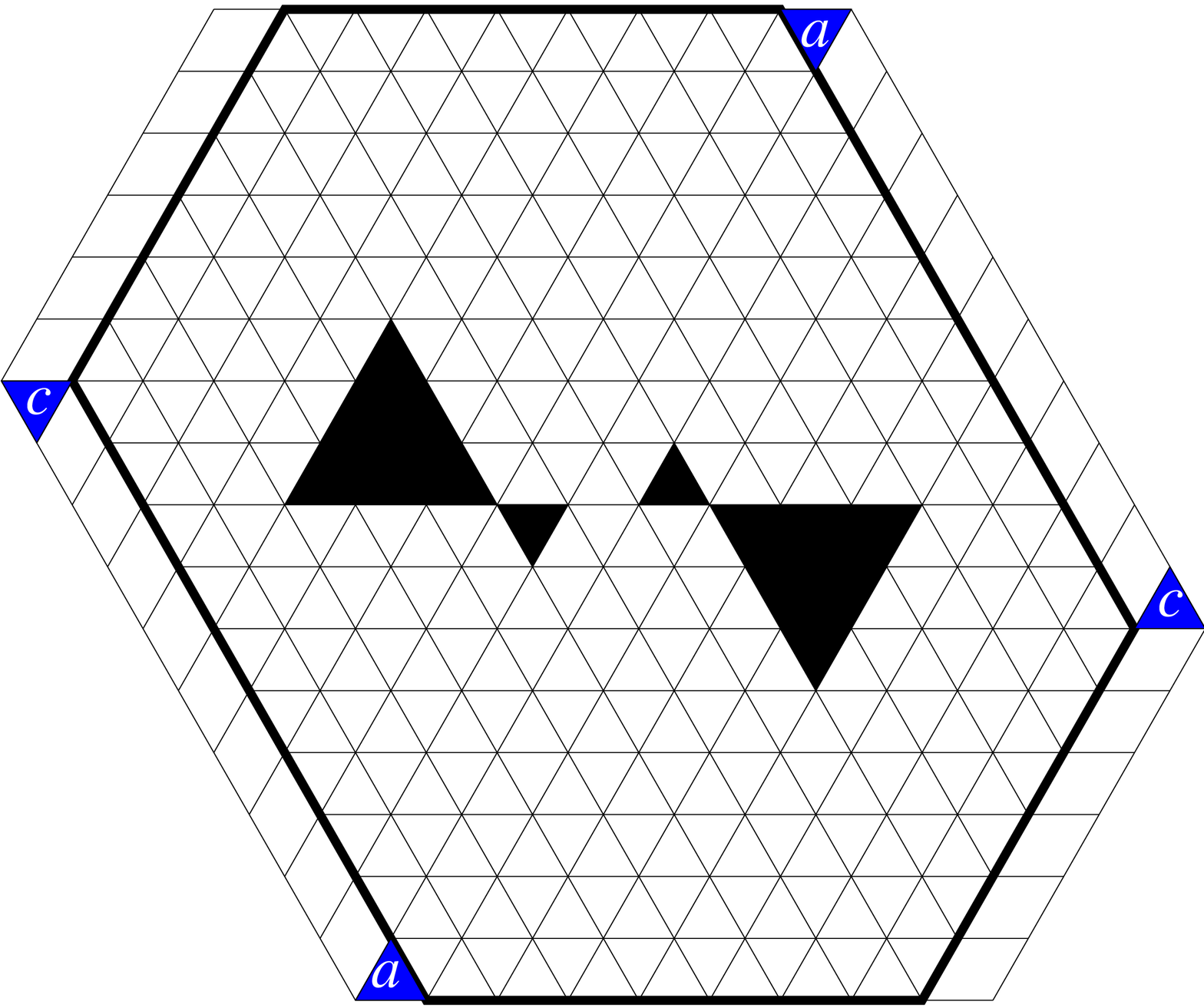}}
\hfill
{\includegraphics[width=0.36\textwidth]{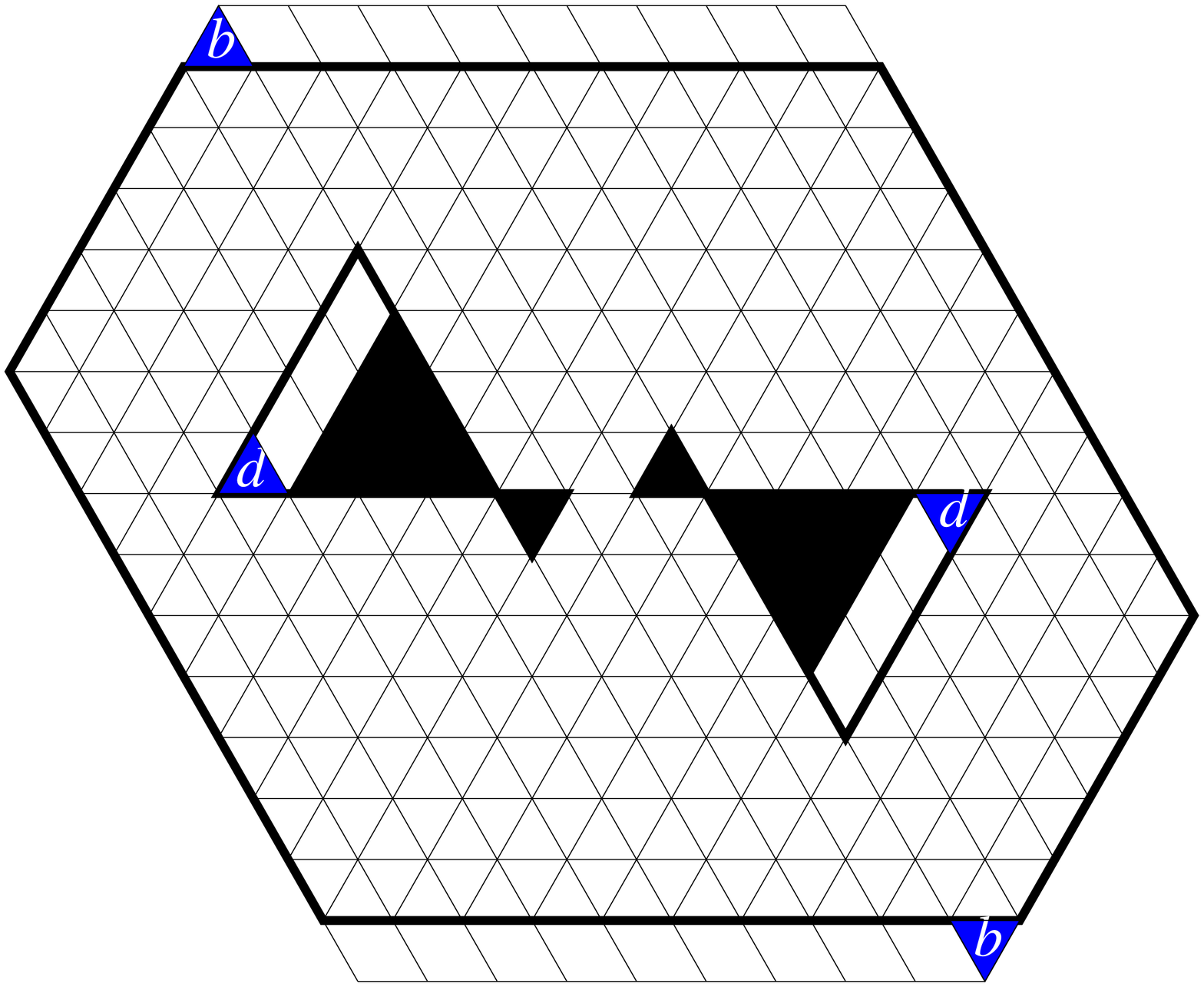}}
\hfill
}
\vskip0.06in
  \centerline{
\hfill
{\includegraphics[width=0.36\textwidth]{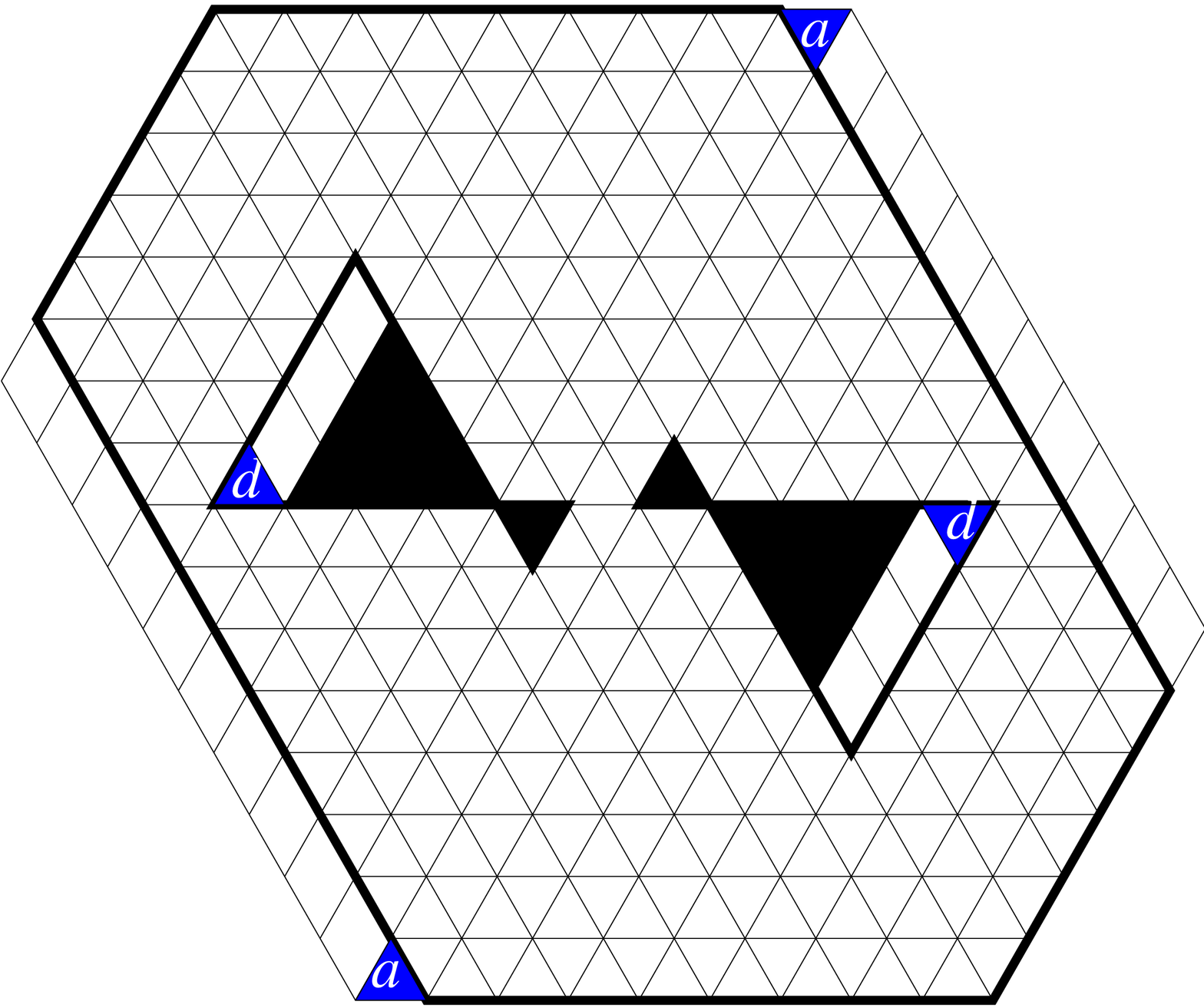}}
\hfill
{\includegraphics[width=0.36\textwidth]{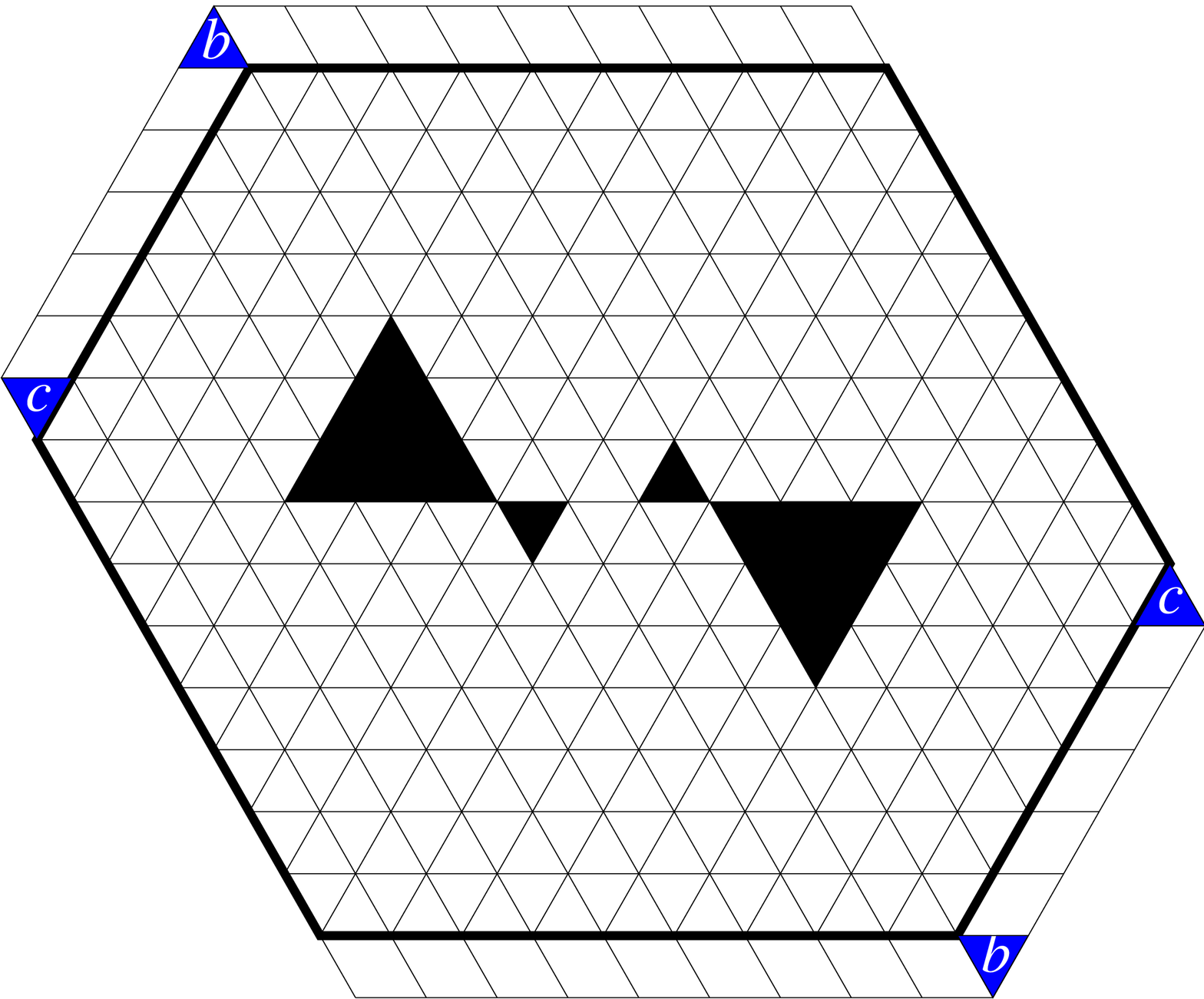}}
\hfill
}
  \caption{\label{fcc} Obtaining the recurrence for the number of centrally symmetric tilings of the regions $FC'_{x,y,z}(a_1,\dotsc,a_k,a_k,\dotsc,a_1)$.}
\vskip-0.35in
\end{figure}

{\it Proof of Theorem \ref{tbb}.} The proof follows the same approach as the proof of Theorem \ref{tba}. 
We consider the planar dual graph $G$ of the region $FC'_{x,y,z}(a_1,\dotsc,a_k,a_k,\dotsc,a_1)$, and apply Theorem \ref{kuovarp} to it, with the vertices $a_1,a_2,b_1,b_2,c_1,c_2,d_1,d_2$ in Theorem \ref{kuovarp} chosen as indicated on the top right in Figure \ref{fcc} (as before, $a_1$ corresponds to the marked unit triangle bottom left, $b_1$, $c_1$, $a_2$, $b_2$, $c_2$ follow in counterclockwise order, and $d_1$ and $d_2$ are next to the removed ferns). 
After removing the forced lozenges, it follows from Figure \ref{fcc} that Theorem \ref{kuovarp} gives

%
\begin{align}
&\!\!
\M_{\odot}(FC'_{x,y,z}(a_1,\dotsc,a_k,a_k,\dotsc,a_1))\M_{\odot}(FC'_{x-2,y-2,z-2}(a_1+1,a_2\dotsc,a_k,a_k,\dotsc,a_2,a_1+1)
\nonumber
\\[5pt]
&
=\M(FC'_{x,y,z-2}(a_1,\dotsc,a_k,a_k,\dotsc,a_1))\M(FC'_{x-2,y-2,z}(a_1+1,a_2,\dotsc,a_k,a_k,\dotsc,a_2,a_1+1))
\nonumber
\\[5pt]
&
+\M_{\odot}(FC'_{x-2,y,z}(a_1,\dotsc,a_k,a_k,\dotsc,a_1))\M_{\odot}(FC'_{x,y-2,z-2}(a_1+1,a_2,\dotsc,a_k,a_k,\dotsc,a_2,a_1+1))
\nonumber
\\[5pt]
&
+
\M_{\odot}(FC'_{x-2,y,z-2}(a_1+1,a_2,\dotsc,a_k,a_k,\dotsc,a_2,a_1+1))\M_{\odot}(FC'_{x,y-2,z}(a_1,\dotsc,a_k,a_k,\dotsc,a_1)),
\label{ecb}
\end{align}

\begin{figure}[h]
  \centerline{
\hfill
{\includegraphics[width=0.33\textwidth]{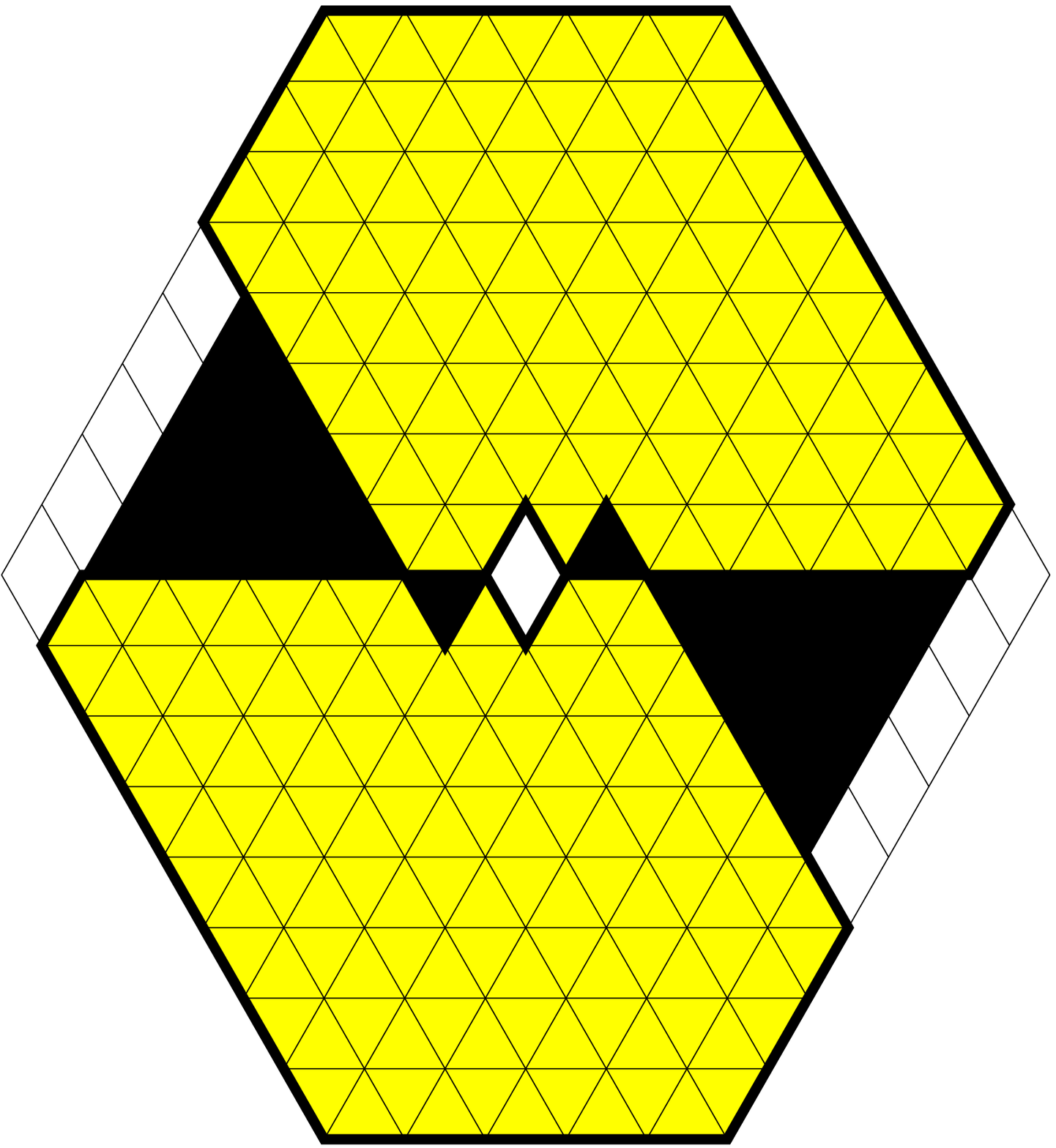}}
\hfill
{\includegraphics[width=0.36\textwidth]{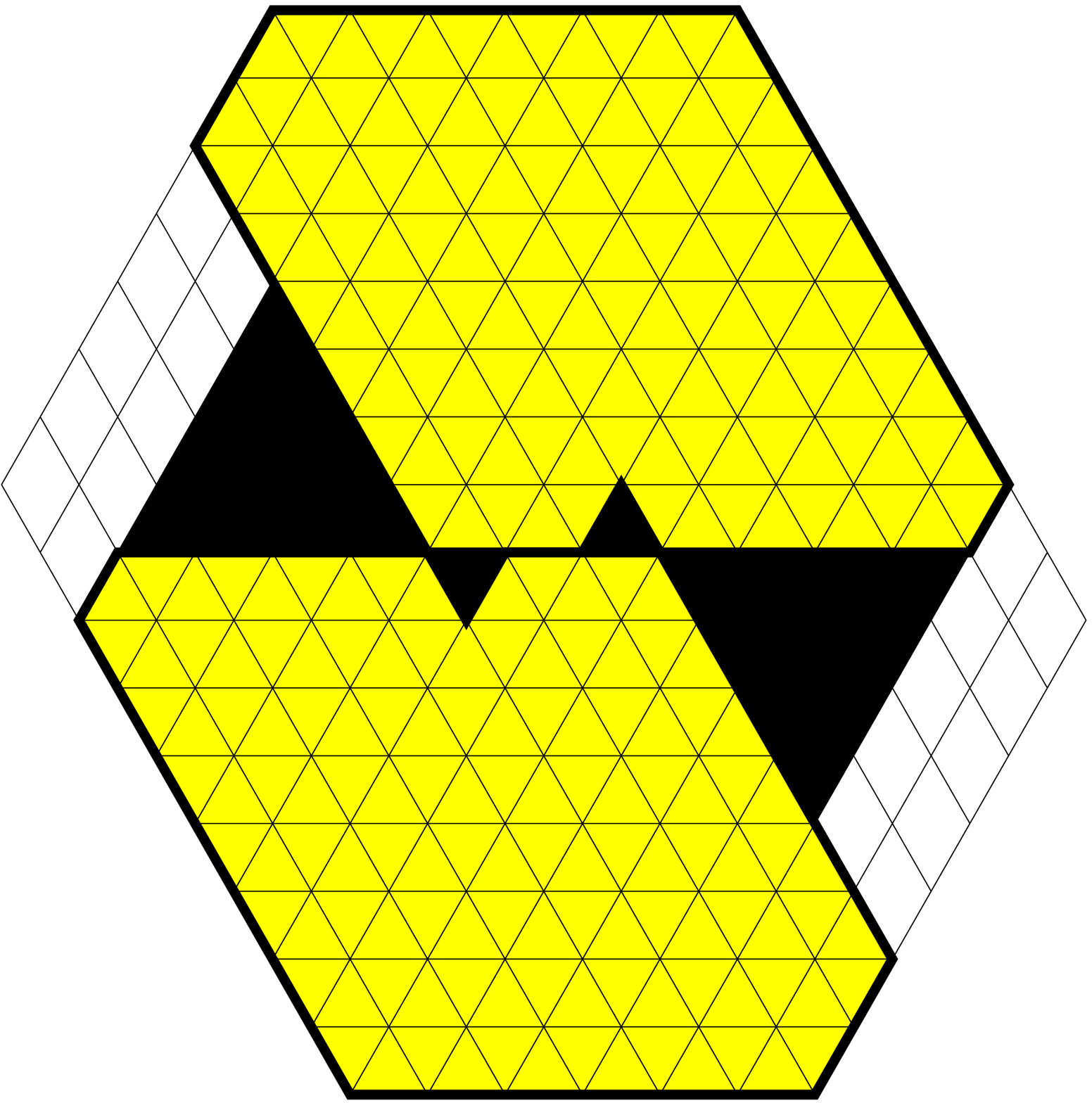}}
\hfill
}
  \caption{\label{fcd} The regions $FC'_{-1,3,3}(4,1,1,4)$ (left) and $FC'_{0,4,2}(4,1,1,4)$ (right).}
\end{figure}

\begin{figure}[h]
  \centerline{
\hfill
{\includegraphics[width=0.42\textwidth]{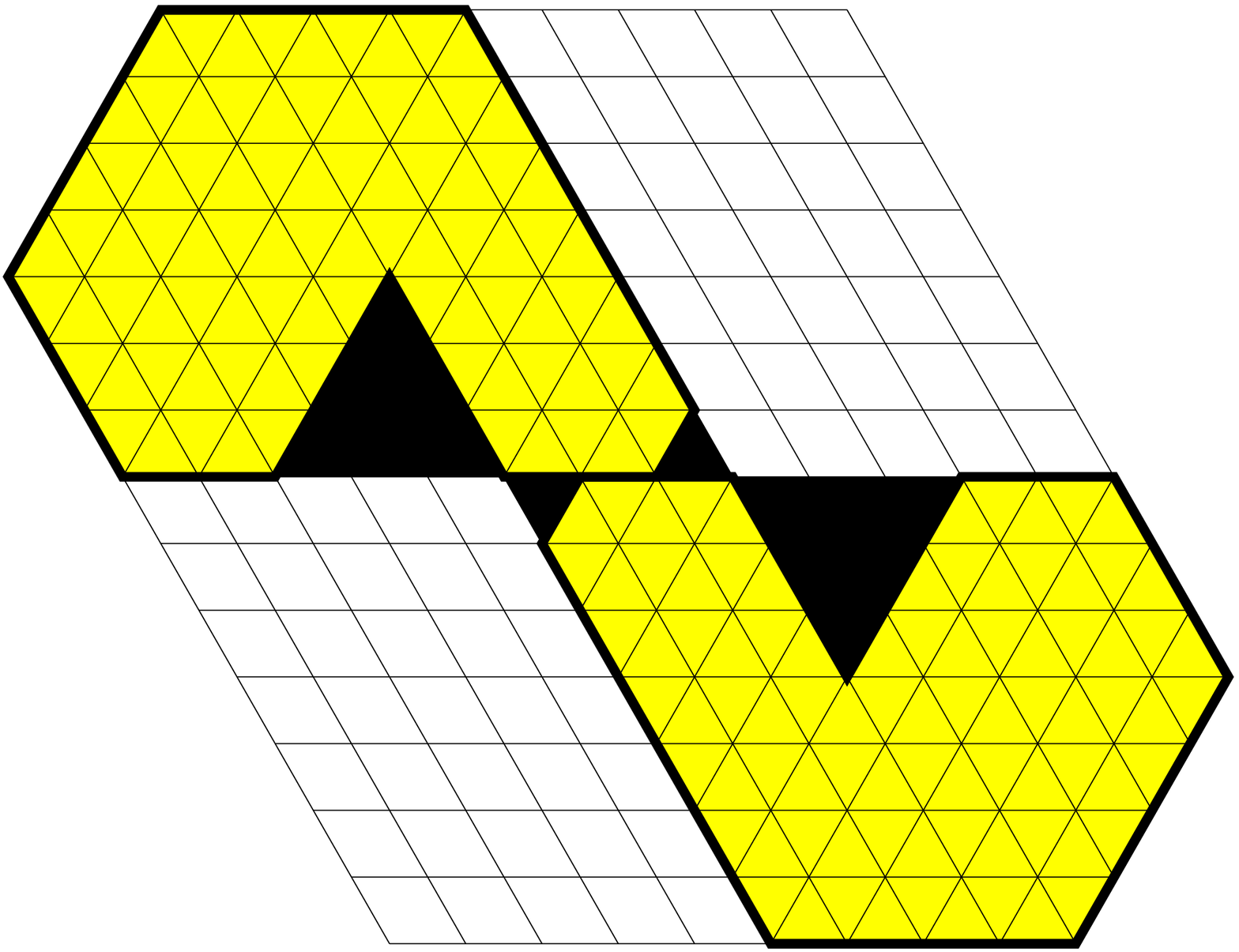}}
\hfill
{\includegraphics[width=0.42\textwidth]{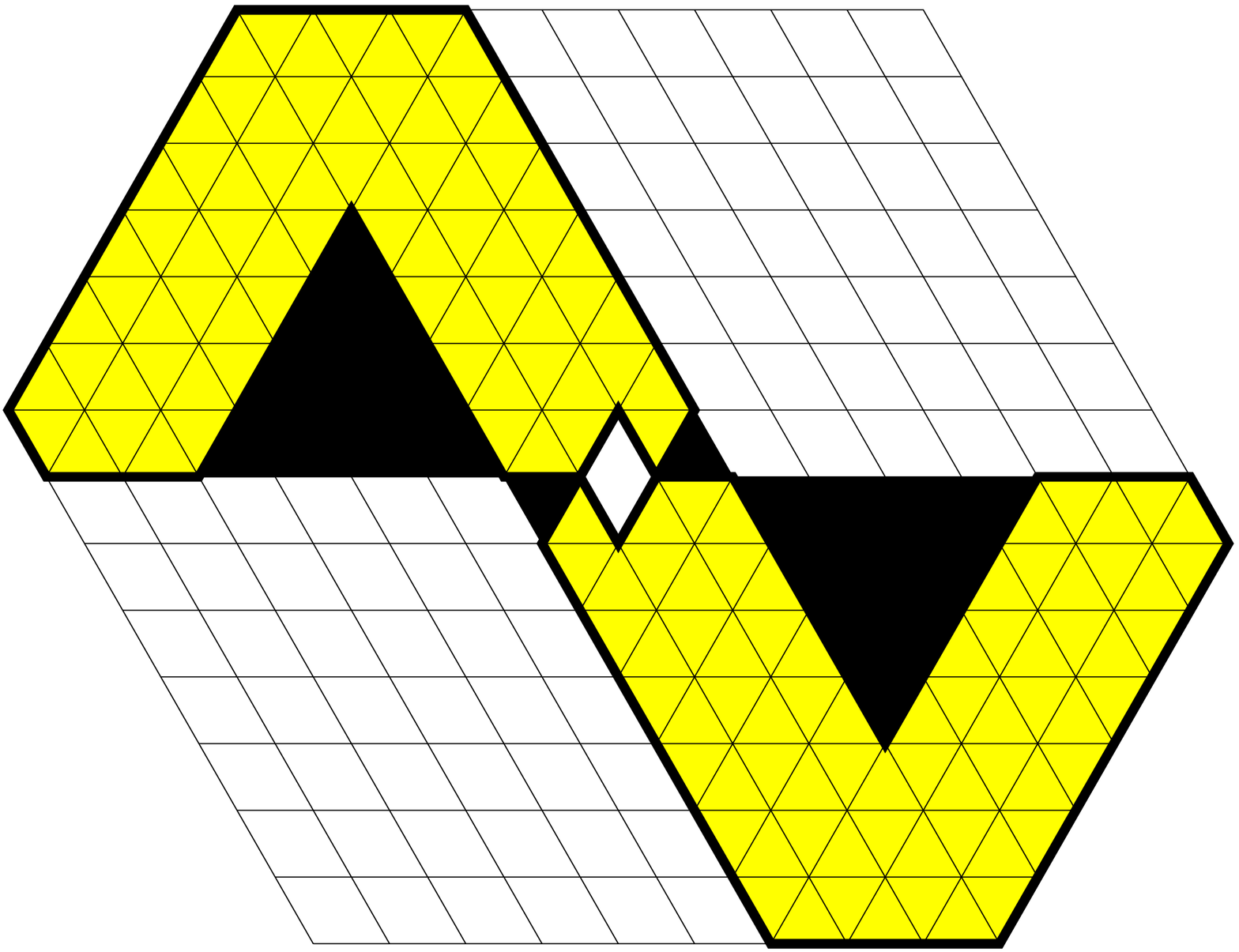}}
\hfill
}
  \caption{\label{fce} The regions $FC'_{4,6,0}(3,1,1,3)$ (left) and $FC'_{3,3,1}(4,1,1,4)$ (right).}
\end{figure}

\parindent0pt
which, when compared to \eqref{eca}, shows that the regions $FC'_{x,y,z}(a_1,\dotsc,a_k,a_k,\dotsc,a_1)$ satisfy precisely the same recurrence as the regions $FC_{x,y,z}(a_1,\dotsc,a_k,a_k,\dotsc,a_1)$.

\parindent15pt
Therefore, by the discution at the beginning of the proof of Theorem \ref{tba}, we can proceed by induction on $x+y+z$. Since $x$, $y$ and $z$ are now only required to have the same parity (and not to be even, as in Theorem \ref{tba}), and since $x\geq-1$ and $y,z\geq0$, the base cases are now $x=-1$, $x=0$, $y=0$, $y=1$, $z=0$ and $z=1$.

Suppose $x=-1$. Then in any tilings of the region $FC'_{-1,y,z}(a_1,\dotsc,a_k,a_k,\dotsc,a_1)$ (illustrated on the left in Figure \ref{fcd} for $y=z=3$, $k=2$, $a_1=4$ and $a_2=1$), the upper shaded region is internally tiled.
Indeed, this follows because the length of its top side is equal to the total length of its base\footnote{ This implies in particular that the central white lozenge in the picture on the left in Figure \ref{fcd} must be part of each tiling of $FC'_{-1,y,z}(a_1,\dotsc,a_k,a_k,\dotsc,a_1)$.}. As in the proof of Theorem \ref{tba}, this, together with formula \eqref{ebacc}, leads to an explicit expression for the $x=-1$ specialization of the left hand side of \eqref{ebc}. One readily verifies that this agrees with the $x=-1$ specialization of the right hand side of \eqref{ebc}, thus verifying \eqref{ebc} in this case.

\parindent15pt
The case $x=0$ is illustrated on the right in Figure \ref{fcd}. We claim that the top shaded region must be internally tiled. Indeed, let $\mathcal{T}$ be a tiling of $FC'_{0,y,z}(a_1,\dotsc,a_k,a_k,\dotsc,a_1)$. Then there is a path of lozenges going up from each of the $a=a_1+\cdots+a_k$ horizontal unit segments on the bottom of this shaded region. These $a$ paths must end on the top side, which has length $a+1$. This leaves precisely one other unit horizontal $s$ segment along the top side. Since $\mathcal{T}$ is centrally symmetric, the path of lozenges that goes down from $s$ must pass between the two removed ferns. This proves our claim. The verification of this base case is finished then in the same fashion as before, using formula \eqref{ebacc}.

The base cases $y=0$ and $y=1$ follow by symmetry from the base cases $x=-1$ and $x=0$, respectively.
The remaining base cases of $z=0$ and $z=1$ are handled in the same fashion. They are illustrated in Figure \ref{fce}. In both of them, it follows by arguments like the ones presented in the previous two paragraphs that the shown shaded regions must be internally tiled, and then explicit product formulas for the corresponding specializations of the numerator and denominator on the left hand side of \eqref{ebc} follow from formula \eqref{ebacc}. These then can be used to verified that the left hand side of \eqref{ebc} agrees with the expression on the right hand side of \eqref{ebc} in these remaining base cases.

As the recurrence \eqref{ecb} is the same as the recurrence \eqref{eca}, the induction step follows in precisely the same way as in the proof of Theorem \ref{tba}. There are only minute differences in the verification that the expression on the right hand side of \eqref{ebc} satisfies the same recurrence, as the expression in Theorem \ref{tbb} is nearly the same as the expression in Theorem \ref{tba}. \epf

\section{Concluding remarks}

In this paper we have proved a simple product formula for the number of centrally symmetric lozenge tilings of a fern-cored hexagon (see Theorem \ref{tba}), thus completing the enumeration of symmetry classes of tilings of these regions. Our formula implies the surprising fact that, when suitably normalized, this number is equal to the square root of the total number of lozenge tilings of the fern-cored hexagon. We also presented a variant when the removed core consists of two disjoint ferns separated by a unit distance (see Theorem \ref{tbb}).

We conjecture that the square root phenomenon mentioned in the previous paragraph holds in more generality, when the removed core is the union of an arbitrary collection of ferns (see Conjectures \ref{tbc} and \ref{tbd}).


It would be interesting to have a bijective proof for the square root phenomenon presented in Conjecture \ref{tbd} (or even for the special case presented in Theorem \ref{tba}).












\end{document}